# A survey of random processes with reinforcement[*]

**Robin Pemantle**

*e-mail:* pemantle@math.upenn.edu

**Abstract:** The models surveyed include generalized Pólya urns, reinforced random walks, interacting urn models, and continuous reinforced processes. Emphasis is on methods and results, with sketches provided of some proofs. Applications are discussed in statistics, biology, economics and a number of other areas.

## Contents

---

[*]This is an original survey paper

## 1. Introduction

In 1988 I wrote a Ph.D. thesis entitled "Random Processes with Reinforcement". The first section was a survey of previous work: it was under ten pages. Twenty years later, the field has grown substantially. In some sense it is still a collection of disjoint techniques. The few difficult open problems that have been solved have not led to broad theoretical advances. On the other hand, some nontrivial mathematics is being put to use in a fairly coherent way by communities of social and biological scientists. Though not full time mathematicians, these scientists are mathematically apt, and continue to draw on what theory there is. I suspect much time is lost, google not withstanding, as they sift through the existing literature and folklore in search of the right shoulders to stand on. My primary motivation for writing this survey is to create universal shoulders: a centralized base of knowledge of the three or four most useful techniques, in a context of applications broad enough to speak to any of half a dozen constituencies of users.

    Such an account should contain several things. It should contain a discussion of the main results and methods, with sufficient sketches of proofs to give a pretty good idea of the mathematics involved[1]. It should contain precise pointers to more detailed statements and proofs, and to various existing versions of the results. It should be historically accurate enough not to insult anyone still living, while providing a modern editorial perspective. In its choice of applications it should winnow out the trivial while not discarding what is simple but useful.

    The resulting survey will not have the mathematical depth of many of the Probability Surveys. There is only one nexus of techniques, namely the stochastic approximation / dynamical system approach, which could be called a theory and which contains its own terminology, constructions, fundamental results, compelling open problems and so forth. There would have been two, but it seems that the multitype branching process approach pioneered by Athreya and Karlin has been taken pretty much to completion by recent work of S. Janson.

    There is one more area that seems fertile if not yet coherent, namely reinforcement in continuous time and space. Continuous reinforcement processes are to reinforced random walks what Brownian motion is to simple random walk, that is to say, there are new layers of complexity. Even excluding the hot new subfield

---

[1] In fact, the heading "PROOF:" in this survey means just such a sketch.

of SLE, which could be considered a negatively reinforced process, there are several other self-interacting diffusions and more general continuous-time processes that open up mathematics of some depth and practical relevance. These are not yet at the mature "surveyable" state, but a section has been devoted to an in-progress glimpse of them.

The organization of the rest of the survey is as follows. Section 2 provides an overview of the basic models, primarily urn models, and corresponding known methods of analysis. Section 3 is devoted to urn models, surveying what is known about some common variants. Section 4 collects applications of these models from a wide variety of disciplines. The focus is on useful application rather than on new mathematics. Section 5 is devoted to reinforced random walks. These are more complicated than urn models and therefore less likely to be taken literally in applications, but have been the source of many of the recognized open problems in reinforcement theory. Section 6 introduces continuous reinforcement processes as well as negative reinforcement. This includes the self-avoiding random walk and its continuous limits, which are well studied in the mathematical physics literature, though not yet thoroughly understood.

## 2. Overview of models and methods

Dozens of processes with reinforcement will be discussed in the remainder of this survey. A difficult organizational issue has been whether to interleave general results and mathematical infrastructure with detailed descriptions of individual processes, or instead whether to lay out the bulk of the mathematics, leaving only some refinements to be discussed along with specific processes and applications. Because of the way research has developed, the existing literature is organized mostly by application; indeed, many existing theoretical results are very much tailored to specific applications and are not easily discussed abstractly. It is, however, possible to describe several distinct approaches to the analysis of reinforcement processes. This section is meant to do so, and to serve as a standalone synopsis of available methodology. Thus, only the most basic urn processes and reinforced random walks will be introduced in this section: just enough to fuel the discussion of mathematical infrastructure. Four main analytical methods are then introduced: exchangeability, branching process embedding, stochastic approximation via martingale methods, and results on perturbed dynamical systems that extend the stochastic approximation results. Prototypical theorems are given in each of these four sections, and pointers are given to later sections where further refinements arise.

### *2.1. Some basic models*

The basic building block for reinforced processes is the urn model[2]. A (single-urn) urn model has an urn containing a number of balls of different types. The set of types may be finite or, in the more general models, countably or uncountably infinite; the types are often taken to be colors, for ease of visualization. The number of balls of each type may be a nonnegative integer or, in the more general models, a nonnegative real number.

---

[2]This is a *de facto* observation, not a definition of reinforced processes.

At each time $n = 1, 2, 3, \ldots$ a ball is drawn from the urn and its type noted. The contents of the urn are then altered, depending on the type that was drawn. In the most straightforward models, the probability of choosing a ball of a given type is equal to the proportion of that type in the urn, but in more general models this may be replaced by a different assumption, perhaps in a way that depends on the time or some aspect of the past, there may be more than one ball drawn, there may be immigration of new types, and so forth.

In this section, the discussion is limited to generalized Pólya urn models, in which a single ball is drawn each time uniformly from the contents of the urn. Sections 3 and 4 review a variety of more general single-urn models. The most general discrete-time models considered in the survey have multiple urns that interact with each other. The simplest among these are mean-field models, in which an urn interacts equally with all other urns, while the more complex have either a spatial structure that governs the interactions or a stochastically evolving interaction structure. Some applications of these more complex models are discussed in Section 4.6. We now define the processes discussed in this section.

Some notation in effect throughout this survey is as follows. Let $(\Omega, \mathcal{F}, \mathbb{P})$ be a probability space on which are defined countable many IID random variables uniform on $[0, 1]$. This is all the randomness we will need. Denote these random variables by $\{U_{nk} : n, k \geq 1\}$ and let $\mathcal{F}_n$ denote the $\sigma$-field $\sigma(U_{mk} : m \leq n)$ that they generate. The variables $\{U_{nk}\}_{k \geq 1}$ are the sources of randomness used to go from step $n-1$ to step $n$ and $\mathcal{F}_n$ is the information up to time $n$. In this section we will need only one uniform random variable $U_n$ at each time $n$, so we let $U_n$ denote $U_{n1}$. A notation that will be used throughout is $\mathbf{1}_A$ to denote the indicator function of the event $A$, that is,

$$\mathbf{1}_A(\omega) := \begin{cases} 1 & \text{if } \omega \in A \\ 0 & \text{if } \omega \notin A \end{cases}.$$

Vectors will be typeset in boldface, with their coordinates denoted by corresponding lightface subscripted variables; for example, a random sequence of $d$-dimensional vectors $\{\mathbf{X}_n : n = 1, 2, \ldots\}$ may be written out as $\mathbf{X}_1 := (X_{11}, \ldots, X_{1d})$ and so forth. Expectations $\mathbb{E}(\cdot)$ always refer to the measure $\mathbb{P}$.

*Pólya's urn*

The original Pólya urn model which first appeared in [EP23; Pól31] has an urn that begins with one red ball and one black ball. At each time step, a ball is chosen at random and put back in the urn along with one extra ball of the color drawn, this process being repeated *ad infinitum*. We construct this recursively: let $R_0 = a$ and $B_0 = b$ for some constants $a, b > 0$; for $n \geq 1$, let $R_{n+1} = R_n + \mathbf{1}_{U_{n+1} \leq X_n}$ and $B_{n+1} = B_n + \mathbf{1}_{U_{n+1} > X_n}$, where $X_n := R_n/(R_n + B_n)$. We interpret $R_n$ as the number of red balls in the urn at time $n$ and $B_n$ as the number of black balls at time $n$. Uniform drawing corresponds to drawing a red ball with probability $X_n$ independent of the past; this probability is generated by our source of randomness via the random variable $U_{n+1}$, with the event $\{U_{n+1} \leq X_n\}$ being the event of drawing a red ball at step $n$.

This model was introduced by Pólya to model, among other things, the spread of infectious disease. The following is the main result concerning this model. The best known proofs, whose origins are not certain [Fre65; BK64], are discussed below.

**Theorem 2.1.** *The random variables $X_n$ converge almost surely to a limit $X$. The distribution of $X$ is $\beta(a,b)$, that is, it has density $Cx^{a-1}(1-x)^{b-1}$ where $C = \dfrac{\Gamma(a+b)}{\Gamma(a)\Gamma(b)}$. In particular, when $a = b = 1$ (the case in [EP23]), the limit variable $X$ is uniform on $[0,1]$.*

The remarkable property of Pólya's urn is that is has a random limit. Those outside of the field of probability often require a lengthy explanation in order to understand this. The phenomenon has been rediscovered by researchers in many fields and given many names such as "lock-in" (chiefly in economic models) and "self organization" (physical models and automata).

*Generalized Pólya urns*

Let us generalize Pólya's urn in several quite natural ways. Take the number of colors to be any integer $k \geq 2$. The number of balls of color $j$ at time $n$ will be denoted $R_{nj}$. Secondly, fix real numbers $\{A_{ij} : 1 \leq i, j \leq k\}$ satisfying $A_{ij} \geq -\delta_{ij}$ where $\delta_{ij}$ is the Kronecker delta function. When a ball of color $i$ is drawn, it is replaced in the urn along with $A_{ij}$ balls of color $j$ for $1 \leq j \leq k$. The reason to allow $A_{ii} \in [-1, 0]$ is that we may think of not replacing (or not entirely replacing) the ball that is drawn. Formally, the evolution of the vector $\mathbf{R}_n$ is defined by letting $\mathbf{X}_n := \mathbf{R}_n / \sum_{j=1}^{k} R_{nj}$ and setting $R_{n+1,j} = R_{nj} + A_{ij}$ for the unique $i$ with $\sum_{t<i} X_{nt} < U_{n+1} \leq \sum_{t \leq i} X_{nt}$. This guarantees that $R_{n+1,j} = R_{nj} + A_{ij}$ for all $j$ with probability $X_{ni}$ for each $i$. A further generalization is to let $\{Y_n\}$ be IID random matrices with mean $A$ and to take $R_{n+1,j} = R_{nj} + (Y_n)_{ij}$ where again $i$ satisfies $\sum_{t<i} X_{nt} < U_{n+1} \leq \sum_{t \leq i} X_{nt}$.

I will use the term **generalized Pólya urn scheme** (GPU) to refer to the model where the reinforcement is $A_{ij}$ and the term **GPU with random increments** when the reinforcement $(Y_n)_{ij}$ involves further randomization. Greater generalizations are possible; see the discussion of time-inhomogeneity in Section 3.2. Various older urn models, such as the Ehrenfest urn model [EE07] can be cast as generalized Pólya urn schemes. The earliest variant I know of was formulated by Bernard Friedman [Fri49]. In Friedman's urn, there are two colors; the color drawn is reinforced by $\alpha > 0$ and the color not drawn is reinforced by $\beta$. This is a GPU with

$$A = \begin{pmatrix} \alpha & \beta \\ \beta & \alpha \end{pmatrix}.$$

Let $X_n$ denote $X_{n1}$, the proportion of red balls (balls of color 1). Friedman analyzed three special cases. Later, David Freedman [Fre65] gave a general analysis of Friedman's urn when $\alpha > \beta > 0$. Freedman's first result is as follows (the paper goes on to find regions of Gaussian and non-Gaussian behavior for $(X_n - \frac{1}{2})$).

**Theorem 2.2** ([Fre65, Corollaries 3.1, 4.1 and 5.1]). *The proportion $X_n$ of red balls converges almost surely to $\frac{1}{2}$.*

What is remarkable about Theorem 2.2 is that the proportion of red balls does *not* have a random limit. It strikes many people as counterintuitive, after coming to grips with Pólya's urn, that reinforcing with, say, 1000 balls of the color drawn and 1 of the opposite color should push the ratio eventually to $\frac{1}{2}$ rather than to a random limit or to $\{0, 1\}$ almost surely. The mystery evaporates rapidly with some back-of-the-napkin computations, as discussed in section 2.4, or with the following observation.

Consider now a generalized Pólya urn with all the $A_{ij}$ strictly positive. The expected number of balls of color $j$ added to the urn at time $n$ given the past is $\sum_i X_{ni} A_{ij}$. By the Perron-Frobenius theory, there is a unique simple eigenvalue whose left unit eigenvector $\pi$ has positive coordinates, so it should not after all be surprising that $\mathbf{X}_n$ converges to $\pi$. The following theorem from to [AK68, Equation (33)] will be proved in Section 2.3.

**Theorem 2.3.** *In a GPU with all $A_{ij} > 0$, the vector $\mathbf{X}_n$ converges almost surely to $\pi$, where $\pi$ is the unique positive left eigenvector of $A$ normalized by $|\pi| := \sum_i \pi_i = 1$.*

*Remark.* When some of the $A_{ij}$ vanish, and in particular when the matrix $A$ has a nontrivial Jordan block for its Perron-Frobenius eigenvalue, then more subtleties arise. We will discuss these in Section 3.1 when we review some results of S. Janson.

*Reinforced random walk*

The first reinforced random walk appearing in the literature was the **edge-reinforced random walk** (ERRW) of [CD87]. This is a stochastic process defined as follows. Let $G$ be a locally finite, connected, undirected graph with vertex set $V$ and edge set $E$. Let $v \sim w$ denote the neighbor relation $\{v, w\} \in E(G)$. Define a stochastic process $X_0, X_1, X_2, \ldots$ taking values in $V(G)$ by the following transition rule. Let $\mathcal{G}_n$ denote the $\sigma$-field $\sigma(X_1, \ldots, X_n)$. Let $X_0 = v$ and for $n \geq 0$, let

$$\mathbb{P}(X_{n+1} = w \,|\, \mathcal{G}_n) = \frac{a_n(w, X_n)}{\sum_{y \sim X_n} a_n(y, X_n)} \tag{2.1}$$

where $a_n(x, y)$ is one plus the number of previous times the edge $\{x, y\}$ has been traversed (in either direction):

$$a_n(x, y) := 1 + \sum_{k=1}^{n-1} \mathbf{1}_{\{X_k, X_{k+1}\} = \{x, y\}}. \tag{2.2}$$

Formally, we may construct such a process by ordering the neighbor set of each vertex $v$ arbitrarily $g_1(v), \ldots, g_{d(v)}(v)$ and taking $X_{n+1} = g_i(X_n)$ if

$$\frac{\sum_{t=1}^{i-1} a_n(g_t(X_n), X_n)}{\sum_{t=1}^{d(X_n)} a_n(g_t(X_n), X_n)} \leq U_n < \frac{\sum_{t=1}^{i} a_n(g_t(X_n), X_n)}{\sum_{t=1}^{d(X_n)} a_n(g_t(X_n), X_n)}. \tag{2.3}$$

In the case that $G$ is a tree, it is not hard to find multi-color Pólya urns embedded in the ERRW. For any fixed vertex $v$, the occupation measures of the edges adjacent to $v$, when sampled at the return times to $v$, form a Pólya urn process, $\{\mathbf{X}_n^{(v)} : n \geq 0\}$. The following lemma from [Pem88a] begins the analysis in Section 5.1 of ERRW on a tree.

**Lemma 2.4.** *The urns $\{\mathbf{X}_n^{(v)}\}_{v \in V(G)}$ are jointly independent.*

The **vertex-reinforced random walk** or VRRW, also due to Diaconis and introduced in [Pem88b], is similarly defined except that the edge weights $a_n(g_t(X_n), X_n)$ in equation (2.3) are replaced by the occupation measure at the destination vertices:

$$a_n(g_t(X_n)) := 1 + \sum_{k=1}^{n} \mathbf{1}_{X_k = g_t(X_n)}. \tag{2.4}$$

For VRRW, for ERRW on a graph with cycles, and for the other variants of reinforced random walk that are defined later, there is no representation directly as a product of Pólya urn processes or even generalized Pólya urn processes, but one may find embedded urn processes that interact nontrivially.

We now turn to the various methods of analyzing these processes. These are ordered from the least to the most generalizable.

### *2.2. Exchangeability*

There are several ways to see that the sequence $\{X_n\}$ in the original Pólya's urn converges almost surely. The prettiest analysis of Pólya's urn is based on the following lemma.

**Lemma 2.5.** *The sequence of colors drawn from Pólya's urn is exchangeable. In other words, letting $C_n = 1$ if $R_n = R_{n-1} + 1$ (a red ball is drawn) and $C_n = 0$ otherwise, then the probability of observing the sequence $(C_1 = \epsilon_1, \ldots, C_n = \epsilon_n)$ depends only on how many zeros and ones there are in the sequence $(\epsilon_1, \ldots, \epsilon_n)$ but not on their order.*

PROOF: Let $\sum_{i=1}^{n} \epsilon_i$ be denoted by $k$. One may simply compute the probabilities:

$$\mathbb{P}(C_1 = \epsilon_1, \ldots, C_n = \epsilon_n) = \frac{\prod_{i=0}^{k-1}(R_0 + i) \prod_{i=0}^{n-k-1}(B_0 + i)}{\prod_{i=0}^{n-1}(R_0 + B_0 + i)}. \tag{2.5}$$

$\square$

It follows by de Finetti's Theorem [Fel71, Section VII.4] that $X_n \to X$ almost surely, and that conditioned on $X = p$, the $\{C_1\}$ are distributed as independent Bernoulli random variables with mean $p$. The distribution of the limiting random variable $X$ stated in theorem 2.1 is then a consequence of the formula (2.5) (see, e.g., [Fel71, VII.4] or [Dur04, Section4.3b]).

The method of exchangeability is neither robust nor widely applicable: the fact that the sequence of draws is exchangeable appears to be a stroke of luck. The method would not merit a separate subsection were it not for two further appearances. The first is in the statistical applications in Section 4.2 below. The second is in ERRW. This process turns out to be **Markov-exchangeable** in the sense of [DF80], which allows an explicit analysis and leads to some interesting open questions, also discussed in Section 5 below.

### 2.3. Embedding

*Embedding in a multitype branching process*

Let $\{\mathbf{Z}(t) := (Z_1(t), \ldots, Z_k(t))\}_{t \geq 0}$ be a branching process in continuous time with $k$ types, and branching mechanism as follows. At all times $t$, each of the $|\mathbf{Z}(t)| := \sum_{i=1}^{k} Z_i(t)$ particles independently branches in the time interval $(t, t+dt]$ with probability $a_i \, dt$. When a particle of type $i$ branches, the collection of particles replacing it may be counted according to type, and the law of this random integer $k$-vector is denoted $\mu_i$. For any $a_1, \ldots, a_k > 0$ and any $\mu_1, \ldots, \mu_k$ with finite mean, such a process is known to exist and has been constructed in, e.g., [INW66; Ath68]. We assume henceforth for nondegeneracy that it is not possible to get from $|\mathbf{Z}(t)| > 0$ to $|\mathbf{Z}(t)| = 0$ and that it is possible to go from $|Z_t| = 1$ to $|Z_t| = n$ for all sufficiently large $n$. We will often also assume that the states form a single irreducible aperiodic class.

Let $0 < \tau_1 < \tau_2 < \cdots$ denote the times of successive branching; our assumptions imply that for all $n$, $\tau_n < \infty = \sup_m \tau_m$. We examine the process $\mathbf{X}_n := \mathbf{Z}(\tau_n)$. The evolution of $\{\mathbf{X}_n\}$ may be described as follows. Let $\mathcal{F}_n = \sigma(\mathbf{X}_1, \ldots, \mathbf{X}_n)$. Then

$$\mathbb{P}(\mathbf{X}_{n+1} = \mathbf{X}_n + \mathbf{v} \,|\, \mathcal{F}_n) = \sum_{i=1}^{k} \frac{a_i X_{ni}}{\sum_{j=1}^{k} a_j X_{nj}} F_i(\mathbf{v} + \mathbf{e}_i), .$$

The quantity $\dfrac{a_i X_{ni}}{\sum_{j=1}^{k} a_j X_{nj}}$ is the probability that the next particle to branch will be of type $i$. When $a_i = 1$ for all $i$, the type of the next particle to branch is distributed proportionally to its representation in the population. Thus, $\{\mathbf{X}_n\}$ is a GPU with random increments. If we further require $F_i$ to be deterministic, namely a point mass at some vector $(A_{i1}, \ldots, A_{ik})$, then we have a classical GPU.

The first people to have exploited this correspondence to prove facts about GPU's were Athreya and Karlin in [AK68]. On the level of strong laws, results about $\mathbf{Z}(t)$ transfer immediately to results about $\mathbf{X}_n = \mathbf{Z}(\tau_n)$. Thus, for example, the fact that $\mathbf{Z}(t)e^{-\lambda_1 t}$ converges almost surely to a random multiple of the Perron-Frobenius eigenvector of the mean matrix $A$ [Ath68, Theorem 1] gives a proof of Theorem 2.3. Distributional results about $\mathbf{Z}(t)$ do not transfer to distributional results about $\mathbf{X}_n$ without some further regularity assumptions; see Section 3.1 for further discussion.

*Embedding via exponentials*

A special case of the above multitype branching construction yields the classical Pólya urn. Each particle independently gives birth at rate 1 to a new particle of the same color (or equivalently, disappears and gives birth to two particles of the original color). This provides yet another means of analysis of the classical Pólya urn, and new generalizations follow. In particular, the collective birth rate of color $i$ may be taken to be a function $f(Z_i)$ depending on the number of particles of color $i$ (but on no other color). Sampling at birth times then yields the dynamic $\mathbf{X}_{n+1} = \mathbf{X}_n + \mathbf{e}_i$ with probability $f(X_{ni})/\sum_{j=1}^k f(X_{nj})$. Herman Rubin was the first to recognize that this dynamic may be de-coupled via the above embedding into independent exponential processes. His observations were published by B. Davis [Dav90] and are discussed in Section 3.2 in connection with a generalized urn model.

To illustrate the versatility of embedding, I include an interesting, if not particularly consequential, application. The so-called OK Corral process is a shootout in which, at time $n$, there are $X_n$ good cowboys and $Y_n$ bad cowboys. Each cowboy is equally likely to land the next successful shot, killing a cowboy on the opposite side. Thus the transition probabilities are $(X_{n+1}, Y_{n+1}) = (X_n - 1, Y_n)$ with probability $Y_n/(X_n + Y_n)$ and $(X_{n+1}, Y_{n+1}) = (X_n, Y_n - 1)$ with probability $X_n/(X_n + Y_n)$. The process stops when $(X_n, Y_n)$ reaches $(0, S)$ or $(S, 0)$ for some integer $S > 0$. Of interest is the distribution of $S$, starting from, say the state $(N, N)$. It turns out (see [KV03]) that the trajectories of the OK Corral process are distributed exactly as time-reversals of the Friedman urn process in which $\alpha = 0$ and $\beta = 1$, that is, a ball is added of the color opposite to the color drawn. The correct scaling of $S$ was known to be $N^{3/4}$ [WM98; Kin99]. By embedding in a branching process, Kingman and Volkov were able to compute the leading term asymptotic for individual probabilities of $S = k$ with $k$ on the order of $N^{3/4}$.

### 2.4. Martingale methods and stochastic approximation

Let $\{\mathbf{X}_n : n \geq 0\}$ be a stochastic process in the euclidean space $\mathbb{R}^n$ and adapted to a filtration $\{\mathcal{F}_n\}$. Suppose that $\mathbf{X}_n$ satisfies

$$\mathbf{X}_{n+1} - \mathbf{X}_n = \frac{1}{n}\left(F(\mathbf{X}_n) + \xi_{n+1} + R_n\right), \tag{2.6}$$

where $F$ is a vector field on $\mathbb{R}^n$, $\mathbb{E}(\xi_{n+1} \,|\, \mathcal{F}_n) = 0$ and the remainder terms $R_n \in \mathcal{F}_n$ go to zero and satisfy $\sum_{n=1}^\infty n^{-1}|R_n| < \infty$ almost surely. Such a process is known as a **stochastic approximation** process after [RM51]; they used this to approximate the root of an unknown function in the setting where evaluation queries may be made but the answers are noisy.

Stochastic approximations arise in urn processes for the following reason. The probability distributions, $Q_n$, governing the color of the next ball chosen are typically defined to depend on the content vector $\mathbf{R}_n$ only via its normalization $\mathbf{X}_n$. If $b$ new balls are added to $N$ existing balls, the resulting increment $\mathbf{X}_{n+1} - \mathbf{X}_n$ is exactly $\frac{b}{b+N}(\mathbf{Y}_n - \mathbf{X}_n)$ where $\mathbf{Y}_n$ is the normalized vector of added balls. Since $b$ is of constant order and $N$ is of order $n$, the mean increment is

$$\mathbb{E}(\mathbf{X}_{n+1} - \mathbf{X}_n \,|\, \mathcal{F}_n) = \frac{1}{n}\left(F(\mathbf{X}_n) + O(n^{-1})\right)$$

where $F(\mathbf{X}_n) = b \cdot \mathbb{E}_{Q_n}(\mathbf{Y}_n - \mathbf{X}_n)$. Defining $\xi_{n+1}$ to be the martingale increment $\mathbf{X}_{n+1} - \mathbb{E}(\mathbf{X}_{n+1} \,|\, \mathcal{F}_n)$ recovers (2.6). Various recent analyses have allowed scaling such as $n^{-\gamma}$ in place of $n^{-1}$ in equation (2.6) for $\frac{1}{2} < \gamma \leq 1$, or more generally, in place of $n^{-1}$, any constants $\gamma_n$ satisfying

$$\sum_n \gamma_n = \infty \tag{2.7}$$

and

$$\sum_n \gamma_n^2 < \infty. \tag{2.8}$$

These more general schemes do not arise in urn and related reinforcement processes, though some of these processes require the slightly greater generality where $\gamma_n$ is a random variable in $\mathcal{F}_n$ with $\gamma_n = \Theta(1/n)$ almost surely. Because a number of available results are not known to hold under (2.7)–(2.8), the term **stochastic approximation** will be reserved for processes satisfying (2.6).

Stochastic approximations arising from urn models with $d$ colors have the property that $\mathbf{X}_n$ lies in the simplex $\Delta^{d-1} := \{\mathbf{x} \in (\mathbb{R}^+)^d : \sum_{i=1}^d x_i = 1\}$. The vector field $F$ maps $\Delta^{d-1}$ to $T\Delta := \{\mathbf{x} \in \mathbb{R}^d : \sum_{i=1}^d x_i = 0\}$. In the two-color case ($d = 2$), the $X_n$ take values in $[0, 1]$ and $F$ is a univariate function on $[0, 1]$. We discuss this case now, then in the next subsection take up the geometric issues arising when $d \geq 3$.

**Lemma 2.6.** *Let the scalar process $\{X_n\}$ satisfy (2.7)–(2.8) and suppose $\mathbb{E}(\xi_{n+1}^2 \,|\, \mathcal{F}_n) \leq K$ for some finite $K$. Suppose $F$ is bounded and $F(x) < -\delta$ for $a_0 < x < b_0$ and some $\delta > 0$. Then for any $[a, b] \subseteq (a_0, b_0)$, with probability 1 the process $\{X_n\}$ visits $[a, b]$ only finitely often. The same holds if $F > \delta$ on $(a_0, b_0)$.*

PROOF: by symmetry we need only consider the case $F < -\delta$ on $(a_0, b_0)$. There is a semi-martingale decomposition $X_n = T_n + Z_n$ where

$$T_n = X_0 + \sum_{k=1}^{n} \gamma_n \left( F(X_{k-1}) + R_{k-1} \right)$$

and

$$Z_n = \sum_{k=1}^{n} \gamma_n \xi_n$$

are respectively the predictable and martingale parts of $X_n$. Square summability of the scaling constants (2.8) implies that $Z_n$ converges almost surely. By assumption, $\sum n^{-1} R_n$ converges almost surely. Thus there is an almost surely finite $N(\omega)$ with

$$|Z_n + R_n - (Z_\infty - R_\infty)| < \frac{1}{2} \min\{a - a_0, b_0 - b\}$$

for all $n \geq N$. No segment of the trajectory of $\{X_{N+k}\}$ can increase by more than $\frac{1}{2} \min\{a - a_0, b_0 - b\}$ while staying inside $[a_0, b_0]$. When $N$ is sufficiently large, the trajectory $\{X_{N+k}\}$ may not jump from $[a, b]$ to the right of $b_0$ nor from the left of $a_0$ to $[a, b]$. The lemma then follows from the observation that for $n > N$, the trajectory if started in $[a, b]$ must exit $[(a + a_0)/2, b]$ to the left and may then never return to $[a, b]$. □

**Corollary 2.7.** *If $F$ is continuous then $X_n$ converges almost surely to the zero set of $F$.*

PROOF: consider the sub-intervals $[a, b]$ of intervals $(a_0, b_0)$ on which $F > \delta$ or $F < -\delta$. Countably many of these cover the complement of the zero set of $F$ and each is almost surely excluded from the limit set of $\{X_n\}$. □

This generalizes a result proved by [HLS80]. They generalized Pólya's urn so that the probability of drawing a red ball was not the proportion $X_n$ of red balls in the urn but $f(X_n)$ for some prescribed $f$. This leads to a stochastic approximation process with $F(x) = f(x) - x$. They also derived convergence results for discontinuous $F$ (the arguments for the continuous case work unless points where $F$ oscillates in sign are dense in an interval) and showed

**Theorem 2.8** ([HLS80, Theorem 4.1]). *Suppose there is a point $p$ and an $\epsilon > 0$ with $F(p) = 0$, $F > 0$ on $(p - \epsilon, p)$ and $F < 0$ on $(p, p + \epsilon)$. Then $\mathbb{P}(X_n \to p) > 0$. Similarly, if $F < 0$ on $(0, \epsilon)$ or $F > 0$ on $(1 - \epsilon, 1)$, then there is a positive probability of convergence to 0 or 1 respectively.*

PROOF, IF $F$ IS CONTINUOUS: Suppose $0 < p < 1$ satisfies the hypotheses of the theorem. By Corollary 2.7, $X_n$ converges to the union of $\{p\}$ and $(p - \epsilon, p + \epsilon)^c$. On the other hand, the semi-martingale decomposition shows that if $X_n$ is in a smaller neighborhood of $p$ and $N$ is sufficiently large, then $\{X_{n+k}\}$ cannot escape $(p - \epsilon, p + \epsilon)$. The cases $p = 0$ and $p = 1$ are similar. □

It is typically possible to find more martingales, special to the problem at hand, that help to prove such things. For the Friedman urn, in the case $\alpha > 3\beta$, it is shown in [Fre65, Theorem 3.1] that the quantity $Y_n := C_n(R_n - B_n)$ is a martingale when $\{C_n\}$ are constants asymptotic to $n^{-\rho}$ for $\rho := (\alpha - \beta)/(\alpha + \beta)$. Similar computations for higher moments show that $\liminf Y_n > 0$, whence $R_n - B_n = \Theta(n^\rho)$.

Much recent effort has been spent obtaining some kind of general hypotheses under which convergence can be shown not to occur at points from which the process is being "pushed away". Intuitively, it is the noise of the process that prevents it from settling down at an unstable zero of $F$, but it is difficult to find the right conditions on the noise and connect them rigorously to destabilization of unstable equilibria. The proper context for a full discussion of this is the next subsection, in which the geometry of vector flows and their stochastic analogues is discussed, but we close here with a one-dimensional result that underlies many of the multi-dimensional results. The result was proved in various forms in [Pem88b; Pem90a].

**Theorem 2.9** (nonconvergence to unstable equilibria). *Suppose $\{X_n\}$ satisfies the stochastic approximation equation (2.6) and that for some $p \in (0,1)$ and $\epsilon > 0$, $\operatorname{sgn} F(x) = \operatorname{sgn}(x - p)$ for all $x \in (p - \epsilon, p + \epsilon)$. Suppose further that $\mathbb{E}(\xi_n^+ \,|\, \mathcal{F}_n)$ and $\mathbb{E}(\xi_n^- \,|\, \mathcal{F}_n)$ are bounded above and below by positive numbers when $X_n \in (p - \epsilon, p + \epsilon)$. Then $\mathbb{P}(\mathbf{X}_n \to p) = 0$.*

PROOF:
**Step 1:** it suffices to show that there is an $\epsilon > 0$ such that for every $n$, $\mathbb{P}(X_k \to p \,|\, \mathcal{F}_n) < 1 - \epsilon$ almost surely. **Proof:** A standard fact is that $\mathbb{P}(X_k \to p \,|\, \mathcal{F}_n) \to 1$ almost surely on the event $\{X_k \to p\}$ (this holds for any event $A$ in place of $\{X_k \to p\}$). In particular, if $\mathbb{P}(X_k \to p) = a > 0$ then for any $\epsilon > 0$ there is some $n$ such that $\mathbb{P}(X_k \to p \,|\, \mathcal{F}_n) > 1 - \epsilon$ on a set of measure at least $a/2$. Thus $\mathbb{P}(X_k \to 0) > 0$ is incompatible with $\mathbb{P}(X_k \to p \,|\, \mathcal{F}_n) < 1 - \epsilon$ almost surely for every $n$.

**Step 2:** with probability $\epsilon$, given $\mathcal{F}_n$, $X_{n+k}$ may wander away from $p$ by $cn^{-1/2}$ due to noise. **Proof:** Let $\tau$ be the exit time of the interval $(p - cn^{-1/2}, p + cn^{-1/2})$. Then $\mathbb{E}(X_\tau - p)^2 \le c^2 n^{-1}$. On the other hand, the quadratic variation of $\{(X_{n \wedge \tau} - p)^2\}$ increases by $\Theta(n^{-2})$ at each step, so on $\{\tau = \infty\}$ is $\Theta(n^{-1})$. If $c$ is small enough, then we see that the event $\{\tau = \infty\}$ must fail at least $\epsilon$ of the time.

**Step 3:** with probability $\epsilon$, $X_{\tau+k}$ may then fail to return to $(p - cn^{-1/2}/2, p + cn^{-1/2}/2)$, due to the drift overcoming the noise. **Proof:** Suppose without loss of generality that $X_\tau < p - cn^{-1/2}$. The quadratic variation of the supermartingale $\{X_{\tau+k}\}$ is $O(\tau^{-1})$, hence $O(n^{-1})$. The probability of such a supermartingale increasing by $cn^{-1/2}/2$ is bounded away from 1. $\square$

As an example, apply this to the urn process in [HLS80], choosing the urn function to be given by $f(x) = 3x^2 - 2x^3$. This corresponds to choosing the color of each draw to be the majority out of three draws sampled with replacement. Here, it may easily be seen that $F < 0$ on $(0, \frac{1}{2})$ and $F > 0$ on $(\frac{1}{2}, 1)$. Verifying the hypotheses on $\xi$, we find that convergence to $\frac{1}{2}$ is impossible, so $S_n \to 0$ or 1 almost surely.

## 2.5. Dynamical systems and their stochastic counterparts

In a vein of research spanning the 1990's and continuing through the present, Benaïm and collaborators have formulated an approach to stochastic approximations based on notions of stability for the approximating ODE. This section describes the dynamical system approach. Much of the material here is taken from the survey [Ben99].

*The dynamical system heuristic*

For processes in any dimension obeying the stochastic approximation equation (2.6) there are two natural heuristics. Sending the noise and remainder terms to zero yields a difference equation $\mathbf{X}_{n+1} - \mathbf{X}_n = n^{-1} F(\mathbf{X}_n)$ and approximating $\sum_{k=1}^n k^{-1}$ by the continuous variable $\log t$ yields the differential equation

$$\frac{d\mathbf{X}}{dt} = F(\mathbf{X}) \,. \tag{2.9}$$

The first heuristic is that trajectories of the stochastic approximation $\{\mathbf{X}_n\}$ should approximate trajectories of the ODE $\{\mathbf{X}(t)\}$. The second is that stable trajectories of the ODE should show up in the stochastic system, but unstable trajectories should not.

A complicating factor in the analysis is the possibility that the trajectories of the ODE are themselves difficult to understand or classify. A standard battery of examples from the dynamical systems literature shows that, once the dimension is greater than one, complicated geometry may arise such as spiraling toward cyclic orbits, orbit chains punctuated by fixed points, and even chaotic trajectories. Successful analysis, therefore, must have several components. First, definitions and results are required in order to understand the forward trajectories of dynamical systems; see the notions of $\omega$-limit sets (forward limit sets) and attractors, below. Next, the notion of trajectory must be generalized to take into account perturbation; see the notions of chain recurrence and chain transitivity below. These topological notions must be further generalized to allow for the kind of perturbation created by stochastic approximation dynamics; see the notion of asymptotic pseudotrajectory below. Finally, with the right definitions in hand, one may prove that a stochastic approximation process $\{\mathbf{X}_n\}$ does in fact behave as an asymptotic pseudotrajectory, and one may establish, under the appropriate hypotheses, versions of the stability heuristic.

It should be noted that an early body of literature exists in which simplifying assumptions preclude flows with the worst geometries. The most common simplifying assumption is that $F = -\nabla V$ for some function $V$, which we think of as a potential. In this case, all trajectories of $\mathbf{X}(t)$ lead "downhill" to the set of local minima of $V$. From the viewpoint of stochastic processes obeying (2.6) that arise in reinforcement models, the assumption $F = -\nabla V$ is quite strong. Recall, however, that the original stochastic approximation processes were designed to locate points such as constrained minima [Lju77; KC78], in which case $F$ is the negative gradient of the objective function. Thus, as pointed out in [BH95; Ben99], much of the early work on stochastic approximation processes focused exclusively on geometrically simple cases such as gradient flow [KC78; BMP90] or attraction to a point [AEK83]. Stochastic approximation processes in the absence of Lyapunov functions can and do follow limit cycles; the earliest natural example I know is found in [Ben97].

*Topological notions*

Although all our flows come from differential equations on real manifolds, many of the key notions are purely topological. A **flow** on a topological space $M$ is a continuous map $(t, x) \mapsto \Phi_t(x)$ from $\mathbb{R} \times M$ to $M$ such that $\Phi_0(x) = x$ and $\Phi_{s+t}(x) = \Phi_t(\Phi_s(x))$ (note that negative times are allowed). The relation to ordinary differential equations is that any bounded Lipschitz vector field $F$ on $\mathbb{R}^n$ has unique integral curves and therefore defines a unique flow $\Phi$ for which $(d/dt)\Phi_t(x) = F(\Phi_t(x))$; we call this the flow associated to $F$. We will assume hereafter that $M$ is compact, our chief example being the $d$-simplex in $\mathbb{R}^{d+1}$. The following constructions and results are due mostly to Bowen and Conley and are taken from Conley's CBMS lecture notes [Con78]. The notions of forward (and backward) limit sets and attractors (and repellers) are old and well known.

For any set $Y \subseteq M$, define the **forward limit set** by

$$\omega(Y) := \bigcap_{t \geq 0} \overline{\bigcup_{s > t} \Phi_s(Y)}. \tag{2.10}$$

When $Y = \{y\}$, this is the set of limit points of the forward trajectory form $y$. Limit sets for sample trajectories will be defined in (2.11) below; a key result will be to relate these to the forward limit sets of the corresponding flow. Reversing time in (2.10), the backward limit set is denoted $\alpha(Y)$.

An **attractor** is a set $A$ that has a neighborhood $U$ such that $\omega(U) = A$. A **repeller** is the time-reversal of this, replacing $\omega(U)$ by $\alpha(U)$. The set $\Lambda_0$ of **rest points** is the set $\{x \in M : \Phi_t(x) = x \text{ for all } t\}$.

Conley then defines the **chain relation** on $M$, denoted $\to$. Say that $x \to y$ if for all $t > 0$ and all open covers U of $M$, there is a sequence $x = z_0, z_1, \ldots, z_{n-1}$, $z_n = y$ of some length $n$ and numbers $t_1, \ldots, t_n \geq t$ such that $\Phi_{t_i}(z_{i-1})$ and $z_i$ are both in some $U \in \mathcal{U}$. In the metric case, this is easier to parse: one must be able to get from $x$ to $y$ by a sequence of arbitrarily long flows separated by arbitrarily small jumps. The **chain recurrent set** $R = R(M, \Phi)$ is defined to be the set $\{x \in M : x \to x\}$. The set $R$ is a compact set containing all **rest points** of the flow (points $x$ such that $\Phi_t(x) = x$ for all $t$), all closures of periodic orbits, and in general all forward and backward limit sets $\omega(y)$ and $\alpha(y)$ of trajectories.

An invariant set $S$ (a union of trajectories) is called **(internally) chain recurrent** if $x \to_S x$ for all $x \in S$, where $\to_S$ denotes the flow restricted to $S$. It is called **(internally) chain transitive** if $x \to_S y$ for all $x, y \in S$. The following equivalence from [Bow75] helps to keep straight the relations between these definitions.

**Proposition 2.10** ([Ben99, Proposition 5.3]). *The following are equivalent conditions on a set $S \subseteq M$.*

1. *$S$ is chain transitive;*
2. *$S$ is chain recurrent and connected;*
3. *$S$ is a closed invariant set and the flow restricted to $S$ has no attractor other than $S$ itself.*

$\square$

**Example 2.1.** Consider the flow on the circle $S^1$ shown on the left-hand side of figure 1. It moves strictly clockwise except at two rest points, $a$ and $b$. Allowing small errors, one need not become stuck at the rest points. The flow is chain recurrent and the only attractor is the whole space. Reversing the flow on the western meridian results in the right-hand figure. Now the point $a$ is a repeller, $b$ is an attractor, the height is a strongly gradient-like function, and the chain recurrent set is $\{a, b\}$.

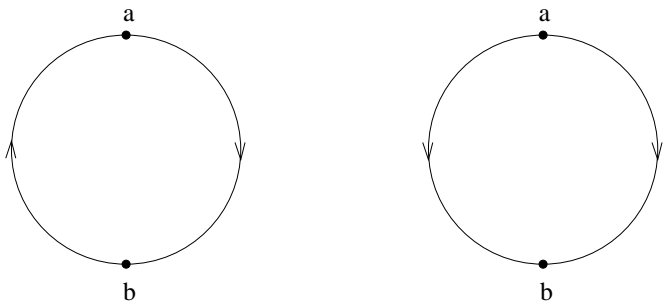

FIG 1. *Two flows on $S^1$*

As we have seen, the geometry is greatly simplified when $F = -\nabla V$. Although this requires differential structure, there is a topological notion that captures the essence. Say that a flow $\{\Phi_t\}$ is **gradient-like** if there is a continuous real function $V : M \to M$ such that $V$ is strictly decreasing along non-constant trajectories. Equation (1) of [Con78, I.5] shows that being gradient-like is strictly weaker than being topologically equivalent to an actual gradient. If in addition, the set $R$ is totally disconnected (hence equal to the set of rest points), then the flow is called **strongly gradient-like**.

Chain recurrence and gradient-like behavior are in some sense the only two possible phenomena. In a gradient-like flow, one can only flow downward. In a chain-recurrent flow, any function weakly decreasing on orbits must in fact be constant on components. Although we will not need the following result, it does help to increase understanding.

**Theorem 2.11** ([Con78, page 17]). *Every flow on a compact space $M$ is uniquely represented as the extension of a chain recurrent flow by a strongly gradient flow. That is, there is a unique subflow (the flow restricted to $R$) which is chain recurrent and for which the quotient flow (collapsing components of $R$ to a point) is strongly gradient-like.* □

*Probabilistic analysis*

An important notion, introduced by Benaïm and Hirsch [BH96], is the **asymptotic pseudotrajectory**. A metric is used in the definition, although it is pointed out in [BLR02, page 13–14] that the property depends only on the topology, not the metric.

**Definition 2.12** (asymptotic pseudotrajectories). *Let $(t,x) \mapsto \Phi_t(x)$ be a flow on a metric space $M$. For a continuous trajectory $X : \mathbb{R}^+ \to M$, let*

$$d_{\Phi,t,T}(X) := \sup_{0 \leq h \leq T} d(X(t+h), \Phi_h(X(t)))$$

*denote the greatest divergence over the time interval $[t, t+T]$ between $X$ and the flow $\Phi$ started from $X(t)$. The trajectory $X$ is an asymptotic pseudotrajectory for $\Phi$ if*

$$\lim_{t \to \infty} d_{\Phi,t,T}(X) = 0$$

*for all $T > 0$.*

This definition is important because it generalizes the "$\to$" relation so that divergence from the flow need not occur at discrete points separated by large times but may occur continuously as long as the divergence remains small over arbitrarily large intervals. This definition also serves as the intermediary between stochastic approximations and chain transitive sets, as shown by the next two results. The first is proved in [Ben99, Proposition 4.4 and Remark 4.5] and the second in [Ben99, Theorem 5.7].

**Theorem 2.13** (stochastic approximations are asymptotic pseudotrajectories). *Let $\{\mathbf{X}_n\}$ be a stochastic approximation process, that is, a process satisfying (2.6), and assume $F$ is Lipschitz. Let $\{\mathbf{X}(t) := \mathbf{X}_n + (t-n)(\mathbf{X}_{n+1} - \mathbf{X}_n)$ for $n \leq t < n+1\}$ linearly interpolate $\mathbf{X}$ at nonintegral times. Assume bounded noise: $|\xi_n| \leq K$. Then $\{\mathbf{X}(t)\}$ is almost surely an asymptotic pseudotrajectory for the flow $\Phi$ of integral curves of $F$.* □

*Remark.* With deterministic step sizes as in (2.6) one may weaken the bounded noise assumption to $L^2$-boundedness: $\mathbb{E}|\xi_n|^2 \leq K$; the stronger assumption is needed only under (2.7)–(2.8). The purpose of the Lipschitz assumption on $F$ is to ensure (along with the standing compactness assumption on $M$) that the flow $\Phi$ is well defined.

The limit set of a trajectory is defined similarly to a forward limit set for a flow. If $X : \mathbb{R}^+ \to M$ is a trajectory, or $X : \mathbb{Z}^+ \to M$ is a discrete time

trajectory, define
$$L(X) := \bigcap_{t \geq 0} \overline{X([t, \infty))}. \tag{2.11}$$

**Theorem 2.14** (asymptotic pseudotrajectories have chain-transitive limits)**.** *The limit set $L(X)$ of any asymptotic pseudotrajectory, $X$, is chain transitive.* □

Combining Theorems 2.13 and 2.14, and drawing on Proposition 2.10 yields a frequently used basic result, appearing first in [Ben93].

**Corollary 2.15.** *Let $X := \{\mathbf{X}_n\}$ be a stochastic approximation process with bounded noise, whose mean vector field $F$ is Lipschitz. Then with probability 1, the limit set $L(X)$ is chain transitive. In view of Proposition 2.10, it is therefore invariant, connected, and contains no proper attractor.* □

Continuing Example 2.1, the right-hand flow has three connected, closed invariant sets $S^1, \{a\}$ and $\{b\}$. The flow restricted to either $\{a\}$ or $\{b\}$ is chain transitive, so either is a possible limit set for $\{\mathbf{X}_n\}$, but the whole set $S^1$ is not chain transitive, thus may not be the limit set of $\{\mathbf{X}_n\}$. We expect to rule out the repeller $\{a\}$ as well, but it is easy to fabricate a stochastic approximation that is rigged to converge to $\{a\}$ with positive probability. Further hypotheses on the noise are required to rule out $\{a\}$ as a limit point. For the left-hand flow, any of the three invariant sets is possible as a limit set.

Examples such as these show that the approximation heuristic, while useful, is somewhat weak without the stability heuristic. Turning to the stability heuristic, one finds better results for convergence than nonconvergence. From [Ben99, Theorem 7.3], we have:

**Theorem 2.16** (convergence to an attractor)**.** *Let $A$ be an attractor for the flow associated to the Lipschitz vector field $F$, the mean vector field for a stochastic approximation $X := \{\mathbf{X}_n\}$. Then either $(i)$ there is a $t$ for which $\{\mathbf{X}_{t+s} : s \geq 0\}$ almost surely avoids some neighborhood of $A$ or $(ii)$ there is a positive probability that $L(X) \subseteq A$ .*

PROOF: A geometric fact requiring no probability is that asymptotic pseudo-trajectories get sucked into attractors. Specifically, let $K$ be a compact neighborhood of the attractor $A$ for which $\omega(K) = A$ (these exist, by definition of an attractor). It is shown in [Ben99, Lemma 6.8] that there are $T, \delta > 0$ such that for any trajectory $X$ starting in $K$, $d_{\Phi,t,T}(X) < \delta$ for all $t$ implies $L(X) \subseteq A$.

Fix such a neighborhood $K$ of $A$ and fix $T, \delta$ as above. By hypothesis, for any $t > 0$ we may find $\mathbf{X}_t \in K$ with positive probability. Theorem 2.13 may be strengthened to yield a $t$ such that

$$\mathbb{P}(d_{\Phi,t,T}(X) < \delta \,|\, \mathcal{F}_t) > 1/2$$

on the event $\mathbf{X}_t \in K$. If $\mathbb{P}(\mathbf{X}_t \in K) = 0$ then conclusion $(i)$ of the theorem is true, while if $\mathbb{P}(\mathbf{X}_t \in K) > 0$, then conclusion $(ii)$ is true. □

For the nonconvergence heuristic, most known results (an exception may be found in [Pem91]) are proved under **linear instability**. This is a stronger hypothesis than topological instability, requiring that at least one eigenvalue of $dF$ have strictly positive real part. An exact formulation may be found in Section 9 of [Ben99]. It is important to that linear instability is defined there for periodic orbits as well as rest points, thus yielding conclusions about nonconvergence to entire orbits, a feature notably lacking in [Pem90a].

**Theorem 2.17** ([Ben99, Theorem 9.1]). *Let $\{\mathbf{X}_n\}$ be a stochastic approximation process on a compact manifold $M$ with bounded noise $||\xi_n|| \leq K$ for all $n$ and $C^2$ vector field $F$. Let $\Gamma$ be a linearly unstable equilibrium or periodic orbit for the flow induced by $F$. Then*

$$\mathbb{P}(\lim_{n \to \infty} d(\mathbf{X}_n, \Gamma) = 0) = 0\,.$$

PROOF: The method of proof is to construct a function $F$ for which $F(\mathbf{X}_n)$ obeys the hypotheses of Theorem 2.9. This relies on known straightening results for stable manifolds and is carried out in [Pem90a] for $\Gamma = \{p\}$ and in [BH95] for general $\Gamma$; see also [Bra98]. □

*Infinite dimensional spaces*

The stochastic approximation processes discussed up to this point obey equation (2.6) which presumes the ambient space $\mathbb{R}^d$. In Section 6.1 we will consider a stochastic approximation on the space $\mathcal{P}(M)$ of probability measures on a compact manifold $M$. The space $\mathcal{P}(M)$ is compact in the weak topology and metrizable, hence the topological definitions of limits, attractors and chain transitive sets are still valid and Theorem 2.14 is still available to force asymptotic pseudotrajectories to have limit sets that are chain transitive. In fact this justifies the space devoted in [Ben99] and its predecessors to establishing results that applied to more than just $\mathbb{R}^d$. The place where new proofs are required is in proving versions of Theorem 2.13 for processes in infinite-dimensional spaces (see Theorem 6.4 below).

*Lyapunov functions*

A **Lyapunov** function for a flow $\Phi$ with respect to the compact invariant set $\Lambda$ is defined to be a continuous function $V : M \to \mathbb{R}$ that is constant on trajectories in $\Lambda$ and strictly decreasing on trajectories not in $\Lambda$. When $\Lambda = \Lambda_0$, the set of rest points, existence of a Lyapunov function is equivalent to the flow being gradient-like. The values $V(\Lambda_0)$ of a Lyapunov function at rest points are called **critical values**. Gradient-like flows are geometrically much better behaved than more general flows, as is shown in [Ben99, Proposition 6.4, and Corollary 6.6]:

**Proposition 2.18** (chain transitive sets when there is a Lyapunov function). *Suppose $V$ is a Lyapunov function for a set $\Lambda$ such that the set of values $V(\Lambda)$* has empty interior. Then every chain transitive set $L$ is contained in $\Lambda$ is a set of constancy for $V$. In particular, if $\Lambda = \Lambda_0$ and intersects the limit set of an asymptotic pseudotrajectory $\{\mathbf{X}(t)\}$ in at most countably many points, then $\mathbf{X}(t)$ must converge to one of these points. $\square$

It follows that the presence of a Lyapunov function for the vector flow associated to $F$ implies convergence of $\{\mathbf{X}_t\}$ to a set of constancy for the Lyapunov function. For example, Corollary 2.7 may be proved by constructing a Lyapunov function with $\Lambda =$ the zero set of $F$. A usual first step in the analysis of a stochastic approximation is therefore to determine whether there is a Lyapunov function. When $F = -\nabla V$ of course $V$ itself is a Lyapunov function with $\Lambda =$ the set of critical points of $V$.

## 3. Urn models: theory

### 3.1. Time-homogeneous generalized Pólya urns

Recall from Section 2.1 the definition of a generalized Pólya urn with reinforcement matrix $A$. We saw in Section 2.3 that the resulting urn process $\{\mathbf{X}_n\}$ may be realized as a multitype branching process $\{\mathbf{Z}(T)\}$ sampled at its jump times $\tau_n$. Already in 1965, for the special case of the Friedman urn with $A := \begin{pmatrix} \alpha & \beta \\ \beta & \alpha \end{pmatrix}$, D. Freedman was able to prove the following limit laws via martingale analysis.

**Theorem 3.1.** *Let $\rho := (\alpha - \beta)/(\alpha + \beta)$. Then*

(i) *If $\rho > 1/2$ then $n^{-\rho}(R_n - B_n)$ converges almost surely to a nontrivial random variable;*
(ii) *If $\rho = 1/2$ then $(n \log n)^{-1/2}(R_n - B_n)$ converges in distribution to a normal with mean zero and variance $(\alpha - \beta)^2$;*
(iii) *If $0 \neq \rho < 1/2$ then $n^{-1/2}(R_n - B_n)$ converges in distribution to a normal with mean zero and variance $(\alpha - \beta)^2/(1 - 2\rho)$.*

Arguments for these results will be given shortly by means of embedding in branching processes. Freedman's original proof of (*iii*) was via moments, estimating each moment by means of an asymptotic recursion; a readable sketch of this argument may be found in [Mah03, Section 6]. The present section summarizes further results that have been obtained via the embedding technique described in Section 2.3. Such an approach rests on an analysis of limit laws in multitype branching processes. These are of independent interest and yet it is interesting to note that such results were not pre-existing. The development of limit laws for multitype branching process was motivated in part by applications to urn processes. In particular, the studies [Ath68] and [Jan04] of multitype limit laws were motivated respectively by the companion paper [AK68] on urn models and by applications to urns in [Jan04; Jan05].

The first thorough study of GPU's via embedding was undertaken by Athreya and Karlin. Although they allow reinforcements to be random, subject to the condition of finite variance, their results depend only on the mean matrix, again denoted $A$. They make an irreducibility assumption, namely that $\exp(tA)$ has positive entries. This streamlines the analysis. While it does not lose too much generality, it probably caused some interesting phenomena in the complementary case to remain hidden for another several decades.

The assumptions imply, by the Perron-Frobenius theory, that the leading eigenvalue of $A$ is real and has multiplicity 1, and that we may write all the eigenvalues as
$$\lambda_1 > \operatorname{Re}\{\lambda_2\} \geq \cdots \geq \operatorname{Re}\{\lambda_d\}.$$

If we do not allow balls to be subtracted and we rule out the trivial case of no reinforcement, then $\lambda_1 > 0$. For any right eigenvector $\xi$ with eigenvalue $\lambda$, the quantity $\xi \cdot \mathbf{Z}(t)e^{-\lambda t}$ is easily seen to be a martingale [AK68, Proposition 1]. When $\operatorname{Re}\{\lambda\} > \lambda_1/2$, this martingale is square integrable, leading to an almost sure limit. This recovers Freedman's first result in two steps. First, taking $\xi = (1,1)$ and $\lambda = \lambda_1 = \alpha + \beta$, we see that $R_n + B_n \sim We^{(\alpha+\beta)t}$ for some random $W > 0$. Secondly, taking $\xi = (1,-1)$ and $\lambda = \alpha - \beta$, we see that $R_n - B_n \sim W'e^{(\alpha-\beta)t}$, with the assumption $\rho > 1/2$ being exactly what is needed square integrability. These two almost sure limit laws imply Freedman's result $(i)$ above.

The analogue of Freedman's result $(iii)$ is that for any eigenvector $\xi$ whose eigenvalue $\lambda$ has $\operatorname{Re}\{\lambda\} < \lambda_1/2$, the quantity $\xi \cdot \mathbf{X}_n/\sqrt{\mathbf{v} \cdot \mathbf{X}_n}$ converges to a normal distribution. The greater generality sheds some light on the reason for the phase transition in the Friedman model at $\rho = 1/2$. For small $\rho$, the mean drift of $R_n - B_n = \mathbf{u} \cdot \mathbf{X}_n$ is swamped by the noise coming from the large number of particles $\mathbf{v} \cdot \mathbf{X}_n = R_n + B_n$. For large $\rho$, early fluctuations in $R_n = B_n$ persist because their mean evolution is of greater magnitude than the noise.

A distributional limit for $\{X_n = \mathbf{Z}(\tau_n)\}$ does not follow automatically from the limit law for $\mathbf{Z}(t)$. A chief contribution of [AK68] is to carry out the necessary estimates to bridge this gap.

**Theorem 3.2** ([AK68, Theorem 3]). *Assume finite variances and irreducibility of the reinforcements. If $\xi$ is a right eigenvector of $A$ whose eigenvalue $\lambda$ satisfies $\operatorname{Re}\{\lambda\} < \lambda_1/2$ then $\xi \cdot \mathbf{X}_n/\sqrt{\mathbf{v} \cdot \mathbf{X}_n}$ converges to a normal distribution.* □

Athreya and Karlin also state that a similar result may be obtained in the "log" case $\operatorname{Re}\{\lambda\} = \lambda_1/2$, extending Freedman's result $(ii)$, but they do not provide details.

At some point, perhaps not until the 1990's, it was noticed that there are interesting cases of GPU's not covered by the analyses of Athreya and Karlin. In particular, the diagonal entries of $A$ may be between $-1$ and $0$, or enough of the off-diagonal entries may vanish that $\exp(tA)$ has some vanishing entries; essentially the only way this can happen is when the urn is **triangular**, meaning that in some ordering of the colors, $A_{ij} = 0$ for $i > j$.

The special case of **balanced** urns, meaning that the row sums of $A$ are constant, is somewhat easier to analyze combinatorially because the total number of balls in the urn increases by a constant each time. Even when the reinforcement is random with mean matrix $A$, the assumption of balance simplifies the analy-sis. Under the assumption of balance and **tenability** (that is, it is not possible for one of the populations to become negative), a number of analyses have been undertaken, including [BP85], [Smy96] and [Mah03]; see also [MS92; MS95] for applications of two-color balanced urns to random recursive trees, and [Mah98] for a tree application of a three-color balanced urn. Exact solutions to two-color balanced urns exhibit involve number theoretic phenomena which are described in [FGP05].

Without the assumption of balance, results on triangular urns date back at least to [DV97]. Their chief results are for two colors, and their method is to analyze the simultaneous functional equations satisfied by the generating functions. Kotz, Mahmoud and Robert [KMR00] concern themselves with removing the balance assumption, attacking the special case $A = \begin{pmatrix} 1 & 0 \\ 1 & 1 \end{pmatrix}$ by combinatorial means. A martingale-based analysis of the cases $A = \begin{pmatrix} 1 & 0 \\ c & 1 \end{pmatrix}$ and $A = \begin{pmatrix} a & 0 \\ 0 & b \end{pmatrix}$ is hidden in [PV99]. The latter case had appeared in various places dating back to [Ros40], the result being as follows.

**Theorem 3.3** (diagonal urn). *Let $a > b > 0$ and consider a GPU with reinforcement matrix*
$$A = \begin{pmatrix} a & 0 \\ 0 & b \end{pmatrix}.$$
*Then $R_n/B_n^\rho$ converges almost surely to a nonzero finite limit, where $\rho := a/b$.*

PROOF: From branching process theory there are variables $W, W'$ with $e^{-at} R_t \to W$ and $e^{-bt} B_t \to W'$. This implies $R_t/B_t^\rho$ converges to the random variable $W/(W')^\rho$, which gives convergence of $R_n/B_n$ to the same quantity. □

Given the piecemeal approaches to GPU's it is fitting that more comprehensive analyses finally emerged. These are due to Janson [Jan04; Jan05]. The first of these is via the embedding approach. The matrix $A$ may be of any finite size, diagonal entries may be as small as $-1$, and the irreducibility assumption is weakened to the largest eigenvalue $\lambda_1$ having multiplicity 1 and being "dominant". This last requirement is removed in [Jan05], which combines the embedding approach with some computations at times $\tau_n$ via generating functions, thus bypassing the need for converting distributional limit theorems in $\mathbf{Z}(t)$ to the stopping times $\tau_n$. The results, given in terms of projections of $A$ onto various subspaces, are somewhat unwieldy to formulate and will not be reproduced here. As far as I can tell, Janson's results do subsume pretty much everything previously known. For example, the logarithmic scaling result appearing in a crude form in [PV99, Theorem 2.3] and elsewhere was proved as Theorem 1.3 (*iv*) of [Jan05]:

**Theorem 3.4.** *Let $R_n$ and $B_n$ be the counts of the two colors of balls in a Friedman urn with $A = \begin{pmatrix} 1 & 0 \\ c & 1 \end{pmatrix}$. Then the quantity $R_n/(cB_n) - \log B_n$ converges* *almost surely to a random finite limit. Equivalently,*

$$\frac{(\log n)^2}{n}\left(B_n - \frac{n}{c\log n} - \frac{n\log\log n}{c(\log n)^2}\right) \tag{3.1}$$

*converges to a random finite limit.* □

To verify the equivalence of the two versions of the conclusion, found respectively in [PV99] and [Jan05], use the deterministic relation $R_n = R_0 + n + (c-1)(B_n - B_0)$ to see that convergence of $R_n/(cB_n) - \log B_n$ is equivalent to

$$\frac{n}{cB_n} - \log B_n = Z + o(1) \tag{3.2}$$

for some finite random $Z$. Also, both versions of the conclusion imply $\log(n/B_n) = \log\log n + \log c + o(1)$ and $\log\log n = \log\log B_n + o(1)$. It follows then that (3.2) is equivalent to

$$\begin{aligned}
B_n &= \frac{n}{c\log B_n + cZ} \\
&= \frac{n}{c\log n}\left(1 + \frac{\log B_n - \log n}{\log n} + \frac{Z + o(1)}{\log n}\right)^{-1} \\
&= \frac{n}{c\log n}\left(1 + \frac{\log(n/B_n)}{\log n} - \frac{Z + o(1)}{\log n}\right) \\
&= \frac{n}{c\log n}\left(1 + \frac{\log\log n}{\log n} - \frac{Z - \log c + o(1)}{\log n}\right)
\end{aligned}$$

which is equivalent to the convergence of (3.1) to the random limit $c^{-1}(Z - \log c)$.

### 3.2. Some variations on the generalized Pólya urn

*Dependence on time*

The **time-dependent** urn is a two-color urn, where only the color drawn is reinforced; the number of reinforcements added at time $n$ is not independent of $n$ but is given by a deterministic sequence of positive real numbers $\{a_n : n = 0, 1, 2, \ldots\}$. This is introduced in [Pem90b] with a story about modeling American primary elections. Denote the contents by $R_n$, $B_n$ and $X_n = R_n/(R_n + B_n)$ as usual. It is easy to see that $X_n$ is a martingale, and the fact that the almost sure limit has no atoms in the open interval $(0,1)$ may be shown via the same three-step nonconvergence argument used to prove Theorem 2.9. The question of atoms among the endpoints $\{0,1\}$ is more delicate. It turns out there is an exact recurrence for the variance of $X_n$, which leads to a characterization of when the almost sure limit is supported on $\{0,1\}$.

**Theorem 3.5** ([Pem90b, Theorem 2]). *Define $\delta_n := a_n/(R_0 + B_0 + \sum_{j=0}^{n-1} a_j)$ to be the ratio of the $n^{th}$ increment to the volume of the urn before the increment is added. Then $\lim_{n\to\infty} X_n = 1$ almost surely if and only if $\sum_{n=1}^{\infty} \delta_n^2 = \infty$.* □

Note that the almost sure convergence of $X_n$ to $\{0,1\}$ is not the same as convergence of $X_n$ to $\{0,1\}$ with positive probability: the latter but not the former happens when $a_n = n$. It is also not the same as almost surely choosing one color only finitely often. No sharp criterion is known for positive probability of $\lim_{n\to\infty} X_n \in \{0,1\}$, but it is known [Pem90b, Theorem 4] that this cannot happen when $\sup_n a_n < \infty$.

*Ordinal dependence*

A related variation adds $a_n$ red balls the $n^{th}$ time a red ball is drawn and $a'_n$ black balls the $n^{th}$ time a black ball is drawn. As is characteristic of such models, a seemingly small change in the definition leads to an different behavior, and to an entirely different method of analysis. One may in fact generalize so that the $n^{th}$ reinforcement of a black ball is of size $a'_n$, not in general equal to $a_n$. The following result appears in the appendix of [Dav90] and is proved by Rubin's exponential embedding.

**Theorem 3.6** (Rubin's Theorem). *Let $S_n := \sum_{k=0}^{n} a_k$ and $S'_n := \sum_{k=0}^{n} a'_n$. Let $G$ denote the event that all but finitely many draws are red, and $G'$ the event that all but finitely many draws are black. Then*

*(i) If $\sum_{n=0}^{\infty} 1/S_n = \infty = \sum_{n=0}^{\infty} 1/S'_n$ then $\mathbb{P}(G) = \mathbb{P}(G') = 0$;*
*(ii) If $\sum_{n=0}^{\infty} 1/S_n = \infty > \sum_{n=0}^{\infty} 1/S'_n$ then $\mathbb{P}(G') = 1$;*
*(iii) If $\sum_{n=0}^{\infty} 1/S_n, \sum_{n=0}^{\infty} 1/S'_n < \infty$ then $\mathbb{P}(G), \mathbb{P}(G') > 0$ and $\mathbb{P}(G) + \mathbb{P}(G') = 1$.*

PROOF: Let $\{Y_n, Y'_n : n = 0, 1, 2, \ldots\}$ be independent exponential with respective means $1/a_n$ and $1/a'_n$. We think of the sequence $Y_1, Y_1 + Y_2, \ldots$ as successive times of an alarm clock. Let $R(t) = \sup\{n : \sum_{k=0}^{n} Y_k \leq t\}$ be the number of alarms up to time $t$, and similarly let $B(t) = \sup\{n : \sum_{k=0}^{n} Y'_k \leq t\}$ be the number of alarms in the primed variables up to time $t$. If $\{\tau_n\}$ are the successive jump times of the pair $(R(t), B(t))$ then $(R(\tau_n), B(\tau_n))$ is a copy of the Davis-Rubin urn process. The theorem follows immediately from this representation, and from the fact that $\sum_{n=0}^{\infty} Y_n$ is finite if and only if its mean is finite (in which case "explosion" occurs) and has no atoms when finite. □

*Altering the draw*

Mahmoud [Mah04] considers an urn model in which each draw consists of $k$ balls rather than just one. There are $k+1$ possible reinforcements depending on how many red balls there are in the sample. This is related to the model of Hill, Lane and Sudderth [HLS80] in which one ball is added each time but the probability it is red is not $X_n$ but $f(X_n)$ for some function $f : [0,1] \to [0,1]$. The end of Section 2.4 introduced the example of majority draw: if three balls are drawn and the majority is reinforced, then $f(x) = x^3 + 3x^2(1-x)$ is the probability that a majority of three will be red when the proportion of reds is $x$. If one samples with replacement in Mahmoud's model and limits the reinforcement to a single ball, then one obtains another special case of the model of Hill, Lane and Sudderth.

A common generalization of these models is to define a family of probability distributions $\{G_x : 0 \leq x \leq 1\}$ on pairs $(Y, Z)$ of nonnegative real numbers, and to reinforce by a fresh draw from $G_x$ when $X_n = x$. If $G_x$ puts mass $f(x)$ on $(1, 0)$ and $1 - f(x)$ on $(0, 1)$, this gives the Hill-Lane-Sudderth urn; an identical model appears in [AEK83]. If $G_x$ gives probability $\binom{k}{j} x^j (1-x)^{k-j}$ to the pair $(\alpha_{1j}, \alpha_{2j})$ for $0 \leq j \leq k$ then this gives Mahmoud's urn with sample size $k$ and reinforcement matrix $\alpha$.

When $G_x$ are all supported on a bounded set, the model fits in the stochastic approximation framework of Section 2.4. For two-color urns, the dimension of the space is 1, and the vector field is a scalar field $F(x) = \mu(x) - x$ where $\mu(x)$ is the mean of $G_x$. As we have already seen, under weak conditions on $F$, the proportion $X_n$ of red balls must converge to a zero of $F$, with points at which the graph of $F$ crosses the $x$-axis in the downward direction (such as the point $1/2$ in a Friedman urn) occurring as the limit with positive probability and points where the graph of $F$ crosses the $x$-axis in an upward direction (such as the point $1/2$ in the majority vote model) occurring as the limit with probability zero.

Suppose $F$ is a continuous function and the graph of $F$ touches the $x$-axis at $(p, 0)$ but does not cross it. The question of whether $X_n \to p$ with positive probability is then more delicate. On one side of $p$, the drift is toward $p$ and on the other side of $p$ the drift is away from $p$. It turns out that convergence can only occur if $X_n$ stays on the side where the drift is toward $p$, and this can only happen if the drift is small enough. A curve tangent to the $x$-axis always yields small enough drift that convergence is possible. The phase transition occurs when the one-sided derivative of $F$ is $-1/2$. More specifically, it is shown in [Pem91] that (*i*) if $0 < F(x) < (p-x)/(2+\epsilon)$ on some neighborhood $(p - \epsilon, p)$ then $X_n \to p$ with positive probability, while (*ii*) if $F(x) > (p-x)/(2-\epsilon)$ on a neighborhood $(p - \epsilon, p)$ and $F(x) > 0$ on a neighborhood $(p, p+\epsilon)$, then $\mathbb{P}(X_n \to p) = 0$. The proof of (*i*) consists of establishing a power law $p - X_n = \Omega(n^{-\alpha})$, precluding $X_n$ ever from exceeding $p$.

The paper [AEK83] introduces the same model with an arbitrary finite number of colors. When the number of colors is $d + 1$, the state vector $\mathbf{X}_n$ lives in the $d$-simplex $\Delta^d := \{(x_1, \ldots, x_{d+1} \in (\mathbb{R}^+)^{d+1} : \sum x_j = 1\}$. Under relatively strong conditions, they prove convergence with probability 1 to a global attractor. A recent variation by Siegmund and Yakir weakens the hypothesis of a global attractor to allow for finitely many non-attracting fixed points on $\partial \Delta^d$ [SY05, Theorem 2.2]. They apply their result to an urn model in which balls are labeled by elements of a finite group: balls are drawn two at a time, and the result of drawing $g$ and $h$ is to place an extra ball of type $g \cdot h$ in the urn. The result is that the contents of the urn converge to the uniform distribution on the subgroup generated by the initial contents.

All of this has been superseded by the stochastic approximation framework of Benaïm *et al.* While convergence to attractors and nonconvergence to repelling sets is now understood, at least in the hyperbolic case (where no eigenvalue of $dF(p)$ has vanishing real part), some questions still remain. In particular, the estimation of deviation probabilities has not yet been carried out. One may ask, for example, how the probability of being at least $\epsilon$ away from a global attractor at time $n$ decreases with $n$, or how fast the probability of being within $\epsilon$ of a repeller at time $n$ decreases with $n$. These questions appear related to quantitative estimates on the proximity to which $\{X_n\}$ shadows the vector flow $\{X(t)\}$ associated to $F$ (cf. the Shadowing Theorem of Benaïm and Hirsch [Ben99, Theorem 8.9]).

## 4. Urn models: applications

In this section, the focus is on modeling rather than theory. Most of the examples contain no significant new mathematical results, but are chosen for inclusion here because they use reinforcement models (mostly urn models) to explain and predict physical or behavioral phenomena or to provide quick and robust algorithms.

### *4.1. Self-organization*

The term **self-organization** is used for systems which, due to micro-level interaction rules, attain a level of coordination across space or time. The term is applied to models from statistical physics, but we are concerned here with self-organization in dynamical models of social networks. Here, self-organization usually connotes a coordination which may be a random limit and is not explicitly programmed into the evolution rules. The Pólya urn is an example of this: the coordination is the approach of $X_n$ to a limit; the limit is random and its sample values are not inherent in the reinforcement rule.

*Market share*

One very broad application of Pólya-like urn models is as a simplified but plausible micro-level mechanism to explain the so-called "lock-in" phenomenon in industrial or consumer behavior. The questions are why one technology is chosen over another (think of the VHS versus Betamax standard for videotape), why the locations of industrial sites exhibit clustering behavior, and so forth. In a series of articles in the 1980's, Stanford economist W. Brian Arthur proposed urn models for this type of social or industrial process, matching data to the predictions of some of the models. Arthur used only very simple urn models, most of which were not new, but his conclusions evidently resonated with the economics community. The stories he associated with the models included the following.
*Random limiting market share:* Suppose two technologies (say Apple versus IBM) are selectively neutral (neither is clearly better) and enter the market at roughly the same time. Suppose that new consumers choose which of the two to buy in proportion to the numbers already possessed by previous consumers. This is the basic Pólya urn model, leading to a random limiting market share: $X_n \to X$. In the case of Apple computers, the sample value of $X$ is between 10% and 15%. This model is discussed at length in [AEK87].

*Random monopoly:* Still assuming no intrinsic advantage, suppose that economies of scale lead to future adoption rates proportional to a power $\alpha > 1$ of present market share. This particular one-dimensional GPU is of the type in Theorem 2.8 (a Hill-Lane-Sudderth urn) with

$$F(x) = \frac{x^\alpha}{x^\alpha + (1-x)^\alpha} - x. \tag{4.1}$$

The graph of $F$ is shaped as in figure 2 below. The equilibrium at $x = 1/2$ is unstable and $X_n$ converges almost surely to 0 or 1. Which of these two occurs depends on chance fluctuations near the beginning of the run. In fact such qualitative behavior persists even if one of the technologies does have an intrinsic advantage, as long as the shape of $F$ remains qualitatively the same. The possibility of an eventual monopoly by an inferior technology is discussed as well in [AEK87] and in the popular account [Art90]. The particular $F$ of (4.1) leads to interesting quantitative questions as to the time the system can spend in disequilibrium, which are discussed in [CL06b; OS05].

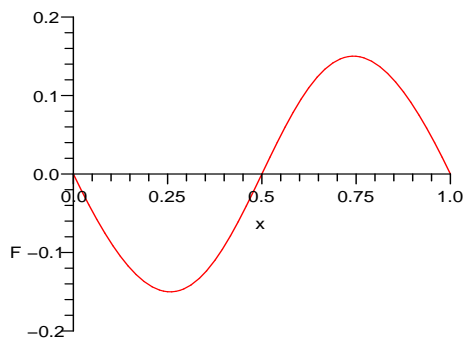

Fig 2. *The urn function $F$ for the power law market share model*

*Neuron polarity*

The mathematics of the following model for neuron growth is mathematically almost identical. The motivating biological question concerns the mechanisms by which apparently identical cells develop into different types. This is poorly understood in many important developmental processes. Khanin and Khanin examine the development of neurons into two types: axon and dendrite. Indistinguishable at first, groups of such cells exhibit periods of growth and retraction until one rapidly elongates to eventually become an axon [KK01, page 1]. They note experimental data suggesting that any neuron has the potential to be either type, and hypotheses that a neuron's length at various stages of growth relative to nearby neurons may influence its development.

They propose an urn model where at each discrete time one of the existing neurons grows by a constant length, $l$, and the others do not grow. The probability of being selected to grow is proportional to the $\alpha$-power of its length, for some parameter $\alpha > 0$. They give rigorous proofs of the long-term behavior in three cases. When $\alpha > 1$, they quote Rubin's Theorem from [Dav90] to show that after a certain random time, only one neuron grows. When $\alpha = 1$, they cite results on the classical Pólya urn from [Fel68] to show that the pairwise length ratios have random finite limits. When $\alpha < 1$, they use embedding methods to show that every pair of lengths has ratio equal to 1 in the limit and to show fluctuations that are Gaussian when $\alpha < 1/2$, Gaussian with a logarithm in the scaling when $\alpha = 1/2$, and differing by a $t^\alpha$ times a random limiting constant when $\alpha \in (1/2, 1)$ (cf. Freedman's results quoted in Section 3.1).

*Preferential attachment*

Another self-organization story has to do with random networks. Models of random networks are used to model the internet, trade, political persuasion and a host of other phenomena. Mathematically, the best studied model is the Erdös-Rényi model where each possible edge is present independently with some probability $p$. For the purposes of many applications, two properties are desirable that do not occur in the Erdös-Rényi model. First, empirical studies show that the distribution of vertex degrees should follow a power law rather than be tightly clustered around its mean. Secondly, there should be local clustering but global connectivity, meaning roughly that as the number of vertices goes to infinity with the average degree constant, the graph-theoretic distance between typical vertices should be small (logarithmic) but the collection of geodesics should have bottlenecks at certain "hub" vertices.

A model, known as the **small-world** model was introduced by Watts and Strogatz [WS98] who were interested in the "six degrees of separation" phenomenon (essentially the empirical fact that the graph of humans and acquaintanceship has local clustering and global connectivity). Their graph is a random perturbation of a nearest neighbor graph. It does exhibit local clustering and global connectivity but not the power-law variation of degrees, and is not easy to work with. A model with the flexibility to fit an arbitrary degree profile was proposed by Chung and Graham and analyzed in [CL03]. This static model is flexible, tractable and provides graphs that match data. Neither this nor the small-world model, however, provides a micro-level explanation of the formation of the graph. A collection of dynamic growth urn models, known as preferential attachment models, the first of which was introduced by Barabási and Albert [BA99], has been developed in order to address this need.

Let a parameter $\alpha \in [0,1]$ be chosen and construct a growing sequence of graphs $\{G_n^\alpha\}$ on the vertex set $\{1,\ldots,n\}$ as follows. Let $G_1$ be the unique graph on one vertex. Given $G_n^\alpha$, let $G_{n+1}^\alpha$ be obtained from $G_n^\alpha$ by adding a single vertex labeled $n+1$ along with a single edge connecting $n+1$ to a random vertex $V_n \in G_n^\alpha$. With probability $\alpha$ the new vertex $V_n$ is chosen uniformly from $\{1,\ldots,n\}$, while with probability $1-\alpha$ the probability $V_n = v$ is taken to be proportional to the degree of $v$.

This procedure always produces a tree. When $\alpha = 1$, this is a well known recursive tree. The other extreme case $\alpha = 0$ may be regarded as pure preferential attachment. A modification is to add some fixed number $m$ of new edges each time, choosing each independently according to the procedure in the case of $m = 1$ and handling collisions among these $m$ new edges by some arbitrary re-sampling scheme. This procedure produces a directed graph that is not, in general, a tree. We denote this random graph by $G_n^{\alpha,m}$.

Preferential attachment models, also known as **rich get richer** models are examples of **scale-free** models [3]. The power laws they exhibit have been fit to data many times, e.g., in figure 1 of [BA99]. Preferential attachment graphs have also been used as the underlying graphs for models of interacting systems. For example, [KKO$^+$05] examines a market pricing model known as the **graphical Fisher model** for price setting. In this model, there is a bipartite graph whose vertices are vendors and buyers. Each buyer buys a unit of goods from the cheapest neighboring vendor, with the vendors trying to set prices as high as possible while still selling all their goods. The emergent prices are entirely a function of the graph structure. In [KKO$^+$05], the graph is taken to be a bipartite version of $G_n^{\alpha,m}$ and the prices are shown to vary only when $m = 1$.

A number of nonrigorous arguments for the degree profile of $G_n^{\alpha,m}$ appear in the literature. For example, in Barabasi and Albert's original paper, the following heuristic argument is given for the case $\alpha = 0$; see also [Mit03]. Consider the vertex $v$ added at time $k$. Let us use an urn model to keep track of its degree. There will be a red ball for each edge incident to $v$ and a black ball for each half of each edge not incident to $v$. The urn begins with $2km$ balls, of which $m$ are red. At each time step a total of $2m$ balls are added. Half of these are always colored black (half-edges incident to $m$ new vertices) while half are colored by choosing from the urn. Let $R_l$ be the number of red balls in the urn at time $l$. Then
$$\mathbb{E}(R_{l+1}|R_l) = R_l \frac{1+m}{2lm} = R_l \frac{l}{2l}$$
and hence
$$\mathbb{E}R_n = m \prod_{l=k}^{n-1}(1 + 1/(2l)) \sim m\sqrt{\frac{n}{k}}.$$

---

[3]see the Wikipedia entry for "scale-free network"

Thus far, the urn analysis is rigorous. The heuristic now proposes that the degree of each ball is exactly the greatest integer below this. Solving for $k$ so that the vertex has degree $d$ at time $n$ gives $k$ as a function of $d$: $k(d) = m^2 n / d^2$. The number of $k$ for which the expected degree is between $d$ and $d+1$ is $\lfloor k(d+1) \rfloor - \lfloor k(d) \rfloor$; this is roughly the derivative with respect to $-d$ of $k(d)$, namely $2m^2 n / d^3$. Thus the fraction of vertices having degree exactly $d$ should be asymptotic to $2m^2 / d^3$.

Chapter 3 of the forthcoming book of Chung and Lu [CL06a] will contain the first rigorous and somewhat comprehensive treatment of preferential attachment schemes (see the discussion in their Section 3.2 of the perils of unjustified heuristics with regard to this model). The only published, rigorous analysis of preferential attachment that I know of is by Bollobás *et al.* [BRST01] and is restricted to the case $\alpha = 0$. Bollobás *et al.* clean up the definition of $G_n^{0,m}$ with regard to the initial conditions and the procedure for resolving collisions. They then prove the following theorem.

**Theorem 4.1** (degrees in the pure preferential attachment graph). *Let*

$$\beta(m,d) := \frac{2m(m+1)}{(m+d)(m+d+1)(m+d+2)}$$

*and let $X_{n,m,d}$ denote the proportion among all $n$ vertices of $G_n^{0,m}$ that have degree $m+d$ (that is, they have in-degree $d$ when edges are directed toward the original vertex). Then both*

$$\inf_{d \leq n^{1/15}} \frac{X_{n,m,d}}{\beta(m,d)}$$

*and*

$$\sup_{d \leq n^{1/15}} \frac{X_{n,m,d}}{\beta(m,d)}$$

*converge to 1 in probability as $n \to \infty$.* □

As $d \to \infty$ with $m$ fixed, $\beta(m,d)$ is asymptotic to $2m^2 d^{-3}$. This agrees, as an asymptotic, with the heuristic for $\alpha = 0$, while providing more information for small $d$. The method of proof is to use Azuma's inequality on the filtration $\sigma(G_n^{0,m} : n = 1, 2, \ldots)$; once this concentration inequality is established, a relatively easy computation finishes the proof by showing convergence of $\mathbb{E} X_{n,m,d}$ to $\beta(m,d)$.

### *4.2. Statistics*

We saw in Theorem 2.1 that the fraction of red balls in a Pólya urn with initial composition $(R(0), B(0))$ converges almost surely and that the limit distribution is $\beta(R(0), B(0))$. Because the sequence of draws is exchangeable, de Finetti's Theorem allows us to interpret the Pólya process as Bayesian observation of a coin with unknown bias, $p$, with a $\beta(R(0), B(0))$ prior on $p$, the probability of flipping "Red" (see the discussion in Section 2.2). Each new flip changes our posterior on $p$, the new posterior after $n$ observations being exactly $\beta(R(n), B(n))$. When $R(0) = B(0) = 1$, the prior is uniform on $[0, 1]$. According to [Fel68, Chapter V, Section 2], Laplace used this model for a tongue-in-cheek estimate that the odds are 1.8 million to one in favor of the sun rising tomorrow; this is based on a record of the sun having risen every day in the modern era (about 5,000 years or 1.8 million days).

*Dirichlet distributions*

The urn representation of the $\beta$ distribution generalizes in the following manner to any number of colors. Consider a $d$-color Pólya urn with initial quantities $R_1(0), \ldots, R_d(0)$. Blackwell and McQueen [BM73, Theorem 1] showed that the limiting distribution is a Dirichlet distribution with parameters $(R_1(0), \ldots, R_d(0))$, where the Dirichlet distribution with parameters $(\alpha_1, \ldots, \alpha_d)$ is defined to be the measure on the $(d-1)$-simplex with density

$$\frac{\Gamma(\alpha_1 + \cdots + \alpha_d)}{\Gamma(\alpha_1) \cdots \Gamma(\alpha_d)} \prod_{j=1}^{d} x_j^{\alpha_j - 1} \, dx_1 \cdots dx_{d-1}. \tag{4.2}$$

The Dirichlet distribution has important statistical properties, some of which we now discuss. Ferguson [Fer73] gives a formula and a discussion of the history. It was long known to Bayesians as the conjugate prior for the parameters of a multinomial distribution (Ferguson refers to [Goo65] for this fact). Thus, for example, the sequence of colors drawn from an urn with initial composition $(1, \ldots, 1)$ are distributed as flips of a $d$-sided coin whose probability vector is drawn from a prior that is uniform on the $(d-1)$-simplex; the posterior after $n$ flips will be a Dirichlet with parameters $(R_1(n), \ldots, R_d(n))$.

Given a finite measure $\alpha$ on a space $S$, the **Dirichlet process** with reference measure $\alpha$ is a random measure $\nu$ on $S$ such that for any disjoint sets $A_1, \ldots, A_d$, the vector of random measures $(\nu(A_1), \ldots, \nu(A_d))$ has a Dirichlet distribution with parameters $(\alpha(A_1), \ldots, \alpha(A_d))$. We denote the law of $\nu$ by $\mathcal{D}(\alpha)$. Because Dirichlet distributions are supported on the unit simplex, the random measure $\nu$ is almost surely a probability measure.

Ferguson [Fer73] suggests using the Dirichlet process as a natural, uninformative prior on the space of probability measures on $S$. Its chief virtue is the ease of computing the posterior: Ferguson shows that after observing independent samples $x_1, \ldots, x_n$ from an unknown measure $\nu$ distributed as $\mathcal{D}(\alpha)$, the posterior for $\nu$ is $\mathcal{D}(\alpha + \sum_{k=1}^{n} \delta(x_k))$, where $\delta(x_k)$ is a point mass at $x_k$. A corollary of this is a beautiful urn representation for $\mathcal{D}(\alpha)$: it is the limiting contents of an $S$-colored Pólya urn with initial "contents" equal to $\alpha$. A second virtue of the Dirichlet prior is that it is weakly dense in the space of probability measures on probability measures on the unit simplex. A drawback is that it is almost surely an atomic measure, meaning that it predicts the eventual occurrence of identical data values. One might prefer a prior supported on the space of continuous measures, although in this regard, the Dirichlet prior is more attractive than its best known predecessor, namely a random distribution function on $[0, 1]$, defined by Dubins and Freedman [DF66], which is almost surely singular-continuous.

The Dirichlet prior and the urn process representing it has been generalized in a number of ways. A random prior on the sequence space $E := \{0, \ldots, k-1\}^\infty$ is defined in [Fer74; MSW92] via an infinite $k$-ary tree of urns. Each urn is a Pólya urn, and the rule for a single update is as follows: sample from the urn at the root; if color $j$ is chosen, put an extra ball of color $j$ in that urn, move to the urn that is the $j^{th}$ child, and repeat this sampling and moving infinitely often. Mapping the space $E$ into any other space $S$ gives a prior on $S$. Taking $k = 2$, $S = [0, 1]$ and the binary map $(x_j) \mapsto \sum x_j 2^{-j}$, one recovers the almost surely singular-continuous prior of [DF66]. Taking $k = 1$, the tree is an infinite ray, and the construction may be used to obtain the Beta-Stacy prior [MSW00].

Another generalization formulates a natural conjugate prior on the the transition matrix of a reversible Markov chain. The edge-reinforced random walk, defined in Section 2.1, is a Markov-exchangeable process (see the last sentence of Section 2.2). This implies that the law of this sequence is a mixture of laws of Markov chains. Given a set of initial weights one the edges, the mixing measure may be explicitly described, as in Theorem 5.1 below. Diaconis and Rolles [DR06] propose this family of such measures, with initial weights as parameters, as priors over reversible Markov transition matrices. Suppose we fix such a prior, coming from initial weights $\{w(e)\}$ and we then observe a single sample $X_0, \ldots, X_n$ of the unknown reversible Markov chain run for time $n$. The posterior distribution will then be another measure from this family, with weights

$$w'(e) := w(e) + \sum_{j=0}^{n-1} \mathbf{1}_{\{X_j, X_{j+1}\}=e} \, .$$

This is exactly analogous to the Ferguson's use of Dirichlet priors for the parameter of an IID sequence and yields, as far as I know, the only computationally feasible Bayesian analysis of an unknown reversible Markov chain.

*The Greenwood-Yule distribution and applications*

Distributions obtained from Pólya urn schemes have been proposed for a variety of applications in which the urn mechanism is plausible at the micro-level. For example, it is proposed in [Jan82] that the number of males born in a family of a specified size $n$ might fit the distribution of a Pólya urn at time $n$ better than a binomial $(n, p)$ if the propensity of having a male was not a constant $p$ but varied according to family. Mackerro and Lawson [ML82] make a similar case (with more convincing data) about the number of days in a given season that are suitable for crop spraying. For more amusing examples, see [Coh76].

Consider a Pólya urn started with $R$ red balls and $n$ black balls and run to time $\alpha n$. The probability that no new balls get added during this time is equal

to
$$\prod_{j=0}^{\alpha n-1} \frac{n+j}{n+R+j}$$

which converges as $n \to \infty$ to $(1+\alpha)^{-R}$. The probability of adding exactly $k$ balls during this time converges as well. To identify the limit, use exchangeability to see that this is $\binom{\alpha n}{k}$ times the probability of choosing zero red balls in $\alpha n - k$ steps and then $k$ red balls in a row. Thus the probability $p_{\alpha n}(k)$ of choosing exactly $k$ red balls is given by

$$p_{\alpha n}(k) = \binom{n}{k} p_{\alpha n - k}(0) \frac{R}{R+(1+\alpha n)-k} \cdots \frac{R+k-1}{R+(1+\alpha)n-1}.$$

The limiting distribution

$$p(k) = (1+\alpha)^{-R} \frac{\prod_{j=0}^{k-1}(R+j)}{(1+\alpha)^k \, k!}$$

is a distribution with very fat tails known as the **Greenwood-Yule distribution** (also, sometimes, the Eggenberger-Pólya distribution). Successive ratios $p(k+1)/p(k)$ are of the form $c\frac{R+k}{k}$, which may be contrasted to the successive ratios $c\frac{R}{k}$ of the Poisson. Thus it is typically used in models where one occurrence may increase the propensity for the next occurrence. It is of historical interest because its use in modeling dependent events precedes the paper [EP23] of Pólya's by several years: the distribution was introduced by Greenwood and Yule [GY20] in order to model numbers of accidents in industrial worksites. More recently it has been proposed as a model for the number of crimes committed by an individual [Gre91], the spontaneous mutation rate in filamentous fungi [BB03] and the number of days in a dry spell [DGVEE05].

It is particularly interesting when the inference process is reversed. The cross-section of the number of particles created in high speed hadronic collisions is known experimentally to have a Greenwood-Yule distribution. This has led physicists [YMN74; Min74] to look for a mechanism responsible for this, perhaps similar to the urn model for Bose-Einstein statistics.

### 4.3. Sequential design

The "two-armed" bandit, whose name seems already to have entered the folklore between 1952 and 1957 [Rob52; BJK62], is a slot machine with two arms. One arm yields a payoff of $1 with probability $p$ and the other arm yields a payoff of $1 with probability $q$. The catch is, you don't know which arm is which, nor do you know $p$ and $q$. The goal is to play so as to maximize your expected return, or limiting average expected return. When $p$ and $q$ are unknown, it is not at all obvious what to do. At the $n^{th}$ step, assuming you have played both arms by then, if you play the arm with the lower historical yield your immediate return is sub-optimal. However, if you always play the arm with the higher historical return, you could miss out forever on a much better action which mis-led you with an initial run of bad luck.

The type of analysis needed to solve the two-armed bandit problem goes by the names of **sequential analysis**, **adaptive control**, or **stochastic** or **optimal** control. Mathematically similar problems occur in statistical hypothesis testing and in the design of clinical trials. The formulation of what is to be optimized, and hence the solution to the problem, will vary with the particular application. In the gambling problem, one wants to maximize the expected return, in the sense of the limiting average (or perhaps the total return in a finite time or infinite time with the future discounted). Determining which of two distributions has a greater mean seems almost identical to the two-armed bandit problem but the objective function is probably some combination of a cost per observation and a reward according to the accuracy of the inference. When designing a clinical trial, say to determine which of two treatments is more effective, there are two competing goals because one is simultaneously gathering data and treating patients. The most data is gathered in a balanced design, where each treatment is tried equally often. But there is an ethical dilemma each time an apparently less effective treatment is prescribed, and the onus is to keep these to a minimum. A survey of both the statistical and ethical problems may be found in [Ros96].

The two-armed bandit problem may be played with asymptotic efficiency. In other words, letting $X_n$ be the payoff at time $n$, there is a strategy such that

$$\lim_{n\to\infty} \frac{1}{n} \sum_{k=1}^{n} X_k = \max\{p, q\}$$

no matter what the values of $p$ and $q$. The first construction I am aware of is due to [Rob52]. A number of papers followed upon that, giving more quantitative solutions in the cases of a finite time horizon [Vog62b; Vog62a], under a finite memory constraint [Rob56; SP65; Sam68], or in a Bayesian framework [Fel62; FvZ70]. One way to formulate an algorithm for asymptotically optimal play is: let $\{\epsilon_n\}$ be a given sequence of real numbers converging to zero; with probability $1 - \epsilon_n$ at time $n$, play whichever arm up to now has the greater average return, and with probability $\epsilon_n$ play the other arm. Such an algorithm is described in [Duf96] and shown to be asymptotically efficient.

In designing a clinical trial, it could be argued that the common good is best served by gathering the most data, since the harm to any finite number of patients who are given the inferior treatment is counterbalanced by the greater efficacy of treatment for all who follow. Block designs, for example alternating between the treatments, were once prevalent but suffer from being predictable by the physician and therefore not double blind.

In 1978, Wei and Durham [WD78] proposed the use of an urn scheme to dictate the sequence of plays in a medical trial. Suppose two treatments have dichotomous outcomes, one succeeding with probability $p$ and the other with probability $q$, both unknown. In Wei and Durham's scheme there is an urn containing at any time two colors of balls, corresponding to the two treatments.

At each time a ball is drawn and replaced, and the corresponding treatment given. If the treatment succeeds, $\alpha$ identical balls and $\beta < \alpha$ balls of the opposite color are added; if the treatment fails, $\alpha$ balls of the opposite color and $\beta$ balls of the same color are added. This is a GPU with random reinforcement and mean reinforcement matrix

$$\begin{pmatrix} p\alpha + (1-p)\beta & (1-p)\alpha + p\beta \\ (1-q)\alpha + q\beta & q\alpha + (1-q)\beta \end{pmatrix}.$$

The unique equilibrium gives nonzero frequencies to both treatments but favors the more effective treatment. It is easy to execute, unpredictable, and comprises between balance and favoring the superior treatment.

If one is relatively more concerned with reducing the number of inferior treatments described, then one seeks something closer to asymptotic efficiency. It is possible to achieve this via an urn scheme as well. Perhaps the simplest way is to reinforce by a constant $\alpha$ if the chosen treatment is effective, but never to reinforce the treatment not chosen. The mean reinforcement matrix for this is simply $\begin{pmatrix} p & 0 \\ 0 & q \end{pmatrix}$. If $p = q$ we have a Pólya urn with a random limit. If $p > q$ we obtain the diagonal urn of Theorem 3.3; the urn population approaches a pure state consisting of only the more effective treatment, with the chance of assigning the inferior treatment at time $n$ being on the order of $n^{-|p-q|/p}$.

Surprisingly, the literature on urn schemes in sequential sampling, as recently as the survey [Dir00] contains no mention of such a scheme. In [LPT04] a stochastic approximation scheme is introduced. Their context is competing investments, and they assume a division of the portfolio into two investments $(X_n, 1 - X_n)$. Let $\{\gamma_n\}$ be a sequence of positive real numbers summing to infinity. Each day, a draw from the urn determines which investment to monitor: the first is monitored with probability $X_n$ and the second with probability $1 - X_n$. If the monitored investment exceeds some threshold, then a fraction $\gamma_n$ of the other investment is transferred into that investment. The respective probabilities for the investments to perform well are unknown and denoted by $p$ and $q$. Defining $T_n$ recursively by $T_n/T_{n+1} = 1 - \gamma_n$, this is a time-dependent Pólya urn process (see Section 3.2) with $a_n = T_{n+1} - T_n$, modified so that the reinforcement only occurs if the chosen investment exceeds the threshold. If $\gamma_n = 1/n$ then $a_n \equiv 1$ and one obtains the diagonal Pólya urn of the preceding paragraph.

When $p \neq q$, the only equilibria are at $X_n = 0$ and $X_n = 1$. The equilibrium at the endpoint 0 is attracting when $p < q$ and repelling when $p > q$, and conversely for the equilibrium at 1. The attractor must be the limit of $\{X_n\}$ with positive probability, but can the repeller be the limit with positive probability? The answer depends on the sequence $\{\gamma_n\}$. It is shown in [LPT04] that for $\gamma_n \sim n^{-\alpha}$, the repeller can be a limit with positive probability when $\alpha < 1$. Indeed, in this case it is easy to see that with positive probability, the attractor is chosen only finitely often. Since we assume $\sum_n \gamma_n = \infty$, this leaves interesting cases near $\gamma_n \approx n^{-1}$. In fact Lamberton, Pagès and Tarrès [LPT04, Corollary 2] show that for $\gamma_n = C/(n+C)$ and $p > q$, the probability of converging to the repeller is zero if and only if $C < 1/p$.

### *4.4. Learning*

A problem of longstanding interest to psychologists is how behavior is learned. Consider a simple model where a subject faces a dichotomous choice: A or B. After choosing, the subject receives a reward. How is future behavior influenced by the reward? Here, the subjects may be animals or humans: in [Her70] pigeons pecked one of two keys and were rewarded with food; in [SP67] the subjects were rats and the reward was pleasant electrical stimulation; in [RE95] the subjects were human and the reward monetary; in [ES54] the subjects were human and success was its own reward. All of these experimenters wished primarily to describe what occurred.

The literature on this sort of learning model is large, but results tend to be mixed, with one model fitting one experiment but not generalizing well. I will, therefore, be content here to describe two popular models and say where they arise. A very basic model is that after a short while, the subject learns which option is best and fixates on that option. According to Herrnstein [Her70, page 243], this does not describe the majority of cases. A hypothesis over 100 years old [Tho98], called the **law of effect**, is that choices will be made with probabilities in proportion to the total reward accumulated when making that choice in the past. Given a (deterministic or stochastic) reward scheme, this then translates into a GPU. In the economic context, the law of effect, also called **the matching law**, is outlined by Roth and Erev [RE95]. They note a resemblance to the evolutionary dynamics formulated by Maynard Smith [MS82], though the models are not the same, and apply their model and some variants to a variety of economic games.

Erev and Roth provide little philosophical justification for the matching law, though their paper has been very influential among evolutionary game theorists. When there are reasons to believe that decision making is operating at a simple level, such models are particularly compelling. In a study of decision making by individuals with brain damage stemming from Huntington's disease, Busemeyer and Stout [BS02] compare a number of plausible models including a Bayesian expected utility model, a stochastic model similar to the Markovian learning models described in the next paragraph, and a Roth-Erev type model. They estimate parameters and test the fit of each model, finding that the Roth-Erev model consistently outperforms the others. See Section 4.6 for more general justifications of this type of model.

A second type of learning model in the psychology literature is a Markovian model with constant step size, which exhibits a stationary distribution rather than convergence to a random limit. Norman [Nor74] reviews several such models, the simplest of which is as follows. A subject repeatedly predicts A or B (in this case, a human predicts whether or not a lamp will flash). The subject's internal state at time $n$ is represented by the probability the subject will choose A, and is denoted $X_n$. The evolution rules contain for parameters, $\theta_1, \ldots, \theta_4 \in (0,1)$. The four possible occurrences are choose A correctly, choose A incorrectly, choose B incorrectly, or choose B correctly, and the new value of $X_{n+1}$ is respectively $X_n + \theta_1(1-X_n)$, $(1-\theta_2)X_n$, $X_n + \theta_3(1-X_n)$ or $(1-\theta_4)X_n$.

Such models were introduced by [ES54; BM55]. The corresponding Markov chain on $[0, 1]$ is amenable to analysis. One interesting result [Nor74, Theorem 3.3] is when $\theta_1 = \theta_4 = \theta$ and $\theta_2 = \theta_3 = 0$. Sending $\theta$ to zero while $n\theta \to t$ gives convergence of $X_{nt}$ to the time-$t$ distribution of a limiting diffusion.

## *4.5. Evolutionary game theory*

Evolutionary game theory is the marriage of the economic concepts of game theory and Nash equilibria with the paradigm of Darwinian evolution originating in biology. A useful reference is [HS98] (replacing the earlier work [HS88]), which has separate introductions for economists and biologists. This subject has exploded in the last several decades, with entire departments and institutes devoted to its study. Naturally, only a very small piece can be discussed here. I will present several applications that reflect the use of urn and reinforcement models, capturing the flavor of this area by giving a vignette rather than a careful history of ideas and methods in evolutionary game theory (and even then, it will take a few pages to arrive at any urn models).

*Economics meets biology*

Applications of evolutionary game theory arise both in economics and biology. This is because each discipline profits considerably from the paradigms of the other, as will now be discussed.

A dominant paradigm in genetics is the stochastic evolution of a genome in a fitness landscape. The **fitness landscape** is a function from genotypes to the real numbers, measuring the adaptive fitness of the corresponding phenotype in the existing environment. A variety of models exist for the change in populations of genotypes based on natural selection with respect to the fitness landscape. Often, randomness is introduced by mechanisms of mutation as well as by stochastic modeling of interactions with the environment. Much of the import of any particular model is in the details of the fitness landscape. Any realistic fitness landscape is hopelessly intractable and different choices of simplifications lead to models illuminating different aspects of evolution.

Game theory enters the biological scene as one type of model for fitness, designed to capture some aspect of the behavior of interacting organisms. Game theoretic models focus on one or two behavioral attributes, usually modeled as expressions of single genes. Different genotypes correspond to different strategies in a single game. Fitness is modeled by the payoff of the given strategy against a mix of other strategies determined by the entire population. Selection acts through increased reproduction as a function of fitness.

In economics, the theory of games and equilibria has been a longstanding dominant paradigm. Interactions between two or more agents are formalized by payoff matrices. Pure and mixed strategies are allowed, but it is generally held that the only strategies that should end up played by rational, informed agents should be **Nash equilibria**[4], that is, strategies that cannot be improved upon given the stochastic mix of strategies in use by the other agents. Two-player games of perfect information are relatively straightforward under assumptions of rationality and perfect information. There is, however, often a distressing lack of correspondence between actual behavior and what is predicted by Nash equilibrium theory.

*Equilibrium selection*

Equilibrium theory can only predict that certain strategies will not be played, leaving open the question of selection among different equilibria. Thus, among the questions that motivated the introduction of evolutionary mechanisms are:

- **equilibrium selection** Which of the equilibria will be played?
- **equilibrium formation** By what route does a population of players come to an equilibrium?
- **equilibrium or not** Will an equilibrium be played at all?

Darwinism enters the economic scene as a means of incorporating bounded information and rationality, explaining equilibrium selection, and modeling games repeated over time and among collections of agents. Assumptions of perfect information and rationality are drastically weakened. Instead, one assumes that individual agents arrive with specific strategies, which they alter only due to data about how well these work (fitness) or to unlikely chance events (mutation). These models make sense in several types of situation. One is when agents are assumed to have low information, for instance in modeling adoption of new technology by consumers, companies, and industries (see the discussion in Section 4.1 of VHS versus Betamax, or Apple versus Mac). Another is when agents are bound by laws, rules or protocols. These, by their nature, must be simple and general[5].

One early application of evolutionary game theory was to explain how players might avoid a Pareto-dominated equilibrium. The ultimate form of this is the Prisoner's dilemma paradox, in which smart people (e.g., game theorists) must choose the only Nash equilibrium, but this is not Pareto-optimal and in fact is dominated by a non-equilibrium play chosen by uneducated people (e.g., mobsters). There are by now many solutions to this dilemma, most commonly involving repeated play. Along the lines of evolutionary game theory, large-scale interactive experiments have been run[6] in which contestants are solicited to submit computer programs that embody various strategies in repeated Prisoner's

---

[4]Many refinements of this notion have been formulated, including subgame-perfect equilibria, coordinated equilibria, etc.

[5]Morals and social norms may be viewed as simple and general principles that may be applied to complex situations. An evolutionary game theoretic approach to explaining these may therefore seem inevitable, and indeed this is the thrust of recent works such as [Sky04; Ale05].

[6]The first was apparently run by Robert Axelrod, a political scientist at the University of Michigan.

Dilemma, and then these are run against each other (in segments of 50 games against each individual opponent) with actual stochastic replicator dynamics to determine which strategies thrive in evolving populations[7].

In the context of more general two-player games, Harsanyi and Selten introduced the concept of the **risk-dominant** equilibrium. This is a notion satisfying certain axioms, among which are naturality not only with respect to game-theoretic equivalences but also the best-reply structure. Consider symmetric $2 \times 2$ games of the form

$$\left( \begin{array}{cc} (a,y) & (0,0) \\ (0,0) & (b,z) \end{array} \right).$$

When $a > b$ and $z > y$ this is a prototypical Nash Bargaining Game. The strategy pair $(1,1)$ is risk-dominant if $ay > bz$. For these games, Pareto-optimality implies risk-dominance, but for other $2 \times 2$ games with multiple equilibria, the risk-dominant equilibrium may not be Pareto-optimal.

Another development in the theory of equilibrium selection, dating back to around 1973, was Selten's **trembling hand**. This is the notion of stochastically perturbing a player's chosen strategy with a small probability $\epsilon$. The idea is that even in an obviously mutually beneficial Nash equilibrium, there is some chance that the opponent will switch to another strategy by mistake (a trembling of the hand), if not through malice or stupidity[8]. A number of notions of equilibria stable under such perturbations arose, depending on the exact model for the $\epsilon$-perturbation, and the way in which $\epsilon \to 0$. An early definition due to J. Maynard Smith was formulated without probability. An **evolutionarily stable strategy** is a strategy such that if it is adopted by a fraction $1 - \epsilon$ of the population, then for sufficiently small $\epsilon$, any other strategy fares worse.

*Replicator dynamics*

One of the earliest and most basic evolutionary game theoretic models is the replicator. There are two versions: the (deterministic) replicator dynamical system and the stochastic replicator. The deterministic replicator assumes a population in which pairs of players with strategy types $1, \ldots, m$ are repeatedly selected at random from a large population, matched against each other in a fixed (generally non-zero-sum) two-player game, and then given a selective advantage in accordance with the outcome of the game. Formally, the model is defined as follows. Fix a two-player (non-zero-sum) game with $m$ strategies for each player such that the payoff to $i$ when playing $i$ against $j$ does not depend on whether the player is Player 1 or Player 2; the matrix of these outcomes is denoted $M$. Let $\mathbf{X}(t)$ denote the normalized population vector, that is, $X_i(t)$ is the proportion of the population at time $t$ that is of type $i$. For any normalized population vector $\mathbf{y}$, the expected outcome for strategy $i$ against a random pick from the population is $E(i, \mathbf{y}) := \sum_{j=1}^{m} M_{i,j} y_j$. Let $E'(i, \mathbf{y}) := E(i, \mathbf{y}) - E_0(\mathbf{y})$ where $E_0(\mathbf{y}) := \sum_{j=1}^{m} E(j, \mathbf{y})$ is the average fitness for the population $\mathbf{y}$; we interpret $\mathbb{E}'(i, \mathbf{y})$ as the selective advantage of type $i$ in population $\mathbf{y}$ and let $\mathbf{E}'(\mathbf{y})$ denote the vector with components $E'(i, \mathbf{y})$. The replicator model is the differential equation

$$\frac{d}{dt}\mathbf{X}(t) = \mathbf{E}'(\mathbf{X}(t)). \tag{4.3}$$

---

[7] One simple strategy that did well in many of these experiments was "Tit for tat": do this time what your opponent did last time.

[8] It is best not to think too much about this when driving past oncoming traffic on a two-lane highway.

The replicator equation was introduced by [TJ78] (though the ideas are present in [MSP73; MS74]) and dubbed "replicator equation" by [SS83]. The books [MS82; HS88] study it extensively.

The notion of evolutionarily stable strategies may be generalized to mixed strategies by means of replicator dynamics. Nash equilibria correspond to rest points for the replicator dynamics. An **evolutionarily stable state** is a population vector that is an attractor for the replicator dynamics (see [HS98, Theorem 7.3.2]).

The presence of the continuous parameter in replicator dynamics indicates that they are a large-population limit. There are a number of discrete systems achieving this limit, but one of the most natural is the **stochastic replicator**. Fix a positive integer $d$ and a $d \times d$ real matrix $M$. We view $M$ as the payoff matrix (for the first player) in a two-player game with $d$ possible strategies, and assume it is normalized to have nonnegative entries. At each integer time $t \geq 0$ there is a population of some size $N(t)$, consisting of individuals whose only attributes are their type, the allowed types being $\{1, \ldots, d\}$. These individuals are represented by an urn with balls of colors $1, \ldots, d$ numbering $N(t)$ altogether. The population at time $t+1$ is determined as follows. Draw $i$ and $j$ at random from the population at time $t$ (with replacement) and return them to the urn along with $M_{ij}$ extra balls of type $i$. The interpretation is that $M_{ij}$ is the fitness of strategy $i$ against strategy $j$ and that the interaction between the two agents causes the representation of type $i$ in the population will change on average by an amount proportional to its fitness against the other strategy it encounters. Repeating this will allow the average growth of type $i$ to be proportional to its average success against all strategies weighted by their representation in the population. One might expect an increase as well of $R_{ji}$ in type $j$, since the interaction has, after all, effects on two agents; in the long run such a term would simply double the rate of change, since an individual will on average be chosen to be Player 1 half the time.

Much of the preceding paragraph is drawn from S. Schreiber's article [Sch01], in which further randomization is allowed ($M$ is the mean matrix for a random increment); as we have seen before, this randomness is not particularly consequential; enough randomness enters through the choice of two individual players. Schreiber also allows $M_{ij} \in [-1, 0]$, which gives his results more general scope than some of their predecessors.

The stochastic replicator is evidently a generalized Pólya urn and its mean ODE is

$$\frac{d\mathbf{Z}(t)}{dt} = \mathbf{Z}(t)^T M \mathbf{Z}(t).$$

One may also consider the normalized population vector $\mathbf{X}(t) := \mathbf{Z}(t)/|\mathbf{Z}(t)|$, where $|\mathbf{Z}(t)|$ is the sum of the components of $\mathbf{Z}(t)$. This evolves, as promised, by a (possibly time-changed) replicator equation

$$\frac{d\mathbf{X}(t)}{dt} = \text{diag}\,(\mathbf{X}(t))\,M\,\mathbf{X}(t) - \mathbf{X}(t)[\mathbf{X}(t)^T\,M\,\mathbf{X}(t)]. \qquad (4.4)$$

In other words, the growth rate of $X_i$ is $X_i(\sum_j M_{ij} X_j - \sum_{rj} X_r M_{ij} X_j)$.

The early study of replicator dynamics concentrated on determining trajectories of the dynamical systems, formulating a notion of stability (such as the **evolutionarily stable strategy** of [MSP73]), and applying these to theoretically interesting biological systems (see especially [MS82]).

The stochastic replicator process fits into the framework of Benaïm *et al.* described in Section 2.5 (except for the possibility of extinction when $M_{ii}$ is allowed to be negative). Schreiber [Sch01, Theorem 2.2] proves a version of Theorem 2.13 for replicator processes, holding on the event of nonextinction. This allows him to derive a version of Corollary 2.15 for replicator process.

**Theorem 4.2** ([Sch01, Corollary 3.2]). *Let $X := \{\mathbf{X}_n\}$ be the normalized population vector for a replicator process with positive expected growth. Then almost surely on the event of nonextinction, the limit set $L(X)$ satisfies the three equivalent properties in Proposition 2.10.* □

It follows from the attractor convergence theorem 2.16 that any attractor in the dynamical system attracts the replicator process with positive probability [BST04, Theorem 7].

Completing the circle ideas, Schreiber has applied his results to a biological model. In [SL96], data is presented showing that three possible color patterns and associated behaviors among the side-blotched lizard *uta stansburiana* have a non-transitive dominance order in terms of success in competing for females[9]. Furthermore, the evolution of population vectors over a six-year period showed a cycle predicted by the dynamical system models of Maynard Smith, which are cited in the paper. Schreiber then applies replicator process urn dynamics. These are the same as in the classic Rock-Paper-Scissors example analyzed in [HS98] and they predict initial cycling followed by convergence to an even mix of all three types in the population.

*Fictitious play*

A quest somewhat related to the problem of explaining equilibrium selection is the problem of finding a mechanism by which a population might evolve toward any equilibrium at all in a game with many strategies. In other words, the emphasis moves from explaining behavior in as Darwinistic a manner as possible to using the idea of natural selection to formulate a coordination algorithm by means of which relatively uninformed agents might adaptively find good (i.e., equilibrium) strategies. Such algorithms are quite important in computer science (internet protocols for use of shared channels, coordination protocols for parallel processing, and so forth).

---

[9]This is evidently the first reported manifestation of this theoretical possibility and I highly recommended reading the brief article to see the details.

In 1951, G. Brown [Bro51] proposed a mechanism known as **fictitious play**. A payoff matrix $M$ is given for a two-player, zero-sum game. Two players play the game repeatedly, with each player choosing at time $n+1$ an action that is optimal if under the assumption that the other player will play according to the past empirical distribution. That is, Player 1 plays $i$ on turn $n+1$ where $i$ is a value of $x$ maximizing the average payoff $n^{-1} \sum_{k=1}^{n} M_{x,y_k}$ and $y_1, \ldots, y_n$ are the previous plays of Player 2; Player 2 plays analogously. Robinson [Rob51] showed that for each player, the empirical distribution of their play converges to an optimal mixed strategy[10].

Fictitious play makes sense for non-zero sum games as well, and for games with more than two players, provided it is specified whether the Bayesian assumption is that each other player independently plays from his empirical distribution or whether the joint play of the other players is from the joint empirical distribution. Robinson's result was extended to non-zero-sum $2 \times 2$ games by [Miy61], but then shown to fail in general by Shapley [Sha64] (a two-player, three-strategy counterexample; see also [Jor93] for a counterexample with dichotomous strategies but three players). There are, however, subclasses of non-zero-sum games for which fictitious play has been shown to converge to Nash equilibria. These include **potential games** [MS96] (every player receives the same payoff), **super-modular games** [MR90] (the payoff matrix is super-modular) and games with interior evolutionarily stable strategies.

Although originally proposed as a computational mechanism, fictitious play became popular behavioral modelers. However, when interpreted as a psychological micro-level mechanism, there are troubling aspects to fictitious play. For a two-player zero-sum game with a unique Nash equilibrium, while the marginals will converge to a saddle point, the plays of the two players may be entirely coordinated, so that actual payoffs may not have the correct long-run average. When there are more than two players, modeling the opponents' future plays as independent picks from empirical marginals seems overly naïve because the empirical joint distribution is known. (The coordination problems that can arise with two players can be thought of in the same way: a failure to model dependence between the opponent's plays one's own plays.) Fudenberg and Kreps [FK93] address these concerns via a greatly generalized framework of optimum response. There chief concern is to give a notion of convergence to Nash equilibrium that precludes the kind of coordination problems mentioned above. In doing so, they take up the notion, due to Harsanyi [Har73], of stochastically perturbed best response, in which each player has independent noise added to the utilities during the computation of the optimum response. They then extend Miyasawa's result on convergence of fictitious play for $2 \times 2$ non-zero-sum games to the setting of stochastic fictitious play, under the assumption of a unique Nash equilibrium [FK93, Proposition 8.1].

---

[10]The same is true with alternating updates, and in fact convergence appears to be faster in that case.

Stochastically perturbed fictitious play fits directly into the stochastic approximation framework. While the stochastic element caused technical difficulties for Fudenberg and Kreps, for whom the available technology was limited to pre-1990 works such as [KC78; Lju77], this same element fits nicely into the framework of Benaïm *et al.* to eliminate unstable trajectories. The groundwork for an analysis in the stochastic approximation framework was laid in [BH99a]. They obtain the usual basic conclusions: the system converges to chain recurrent sets for the associated ODE and attractors attract with positive probability. They give examples of failure to converge, including the stochastic analogue of Jordan's $2 \times 2 \times 2$ counterexample. They then begin to catalogue cases where stochastic fictitious play does converge. Under suitable nondegeneracy assumptions on the noise, they extend [FK93, Proposition 8.1] to allow at most countably many Nash equilibria. Perhaps more interesting is their introduction of a class of two-player $n \times 2$ games they call **generalized coordination games** for which they are able to obtain convergence of stochastic fictitious play. This condition is somewhat restrictive, but in a subsequent work [BH99b], they formulate a simpler and more general condition. Let $F$ denote the vector field of the stochastic approximation process associated with stochastically perturbed fictitious play for a given $m$-player (non-zero-sum) game. Say that $F$ is **cooperative** if $\partial F_i / \partial x_j \geq 0$ for every $i \neq j$. For example, it turns out that the vector field for any generalized coordination game is cooperative. Under a number of technical assumptions, they prove the following result for any cooperative stochastic approximation. Note though, that this is proved for stochastic approximations with constant step size $\epsilon$, as $\epsilon \to 0$; this is in keeping with the prevailing economic formulations of perturbed equilibria, but in contrast to the usual stochastic approximation framework.

**Theorem 4.3** ([BH99b, Theorem 1.5]). *If $F$ is cooperative then as $\epsilon \to 0$, the empirical measure of the stochastic approximation process converges in probability to the set of equilibria of the vector field $F$. If in addition either $F$ is real analytic or has only finitely many stable equilibria, then the empirical distribution converges to an asymptotically stable equilibrium.* □

*Remark.* This result requires constant step size (2.6) but is conjectured to hold under (2.7)–(2.8); see [Ben00, Conjecture 2.3]. The difficulty is that the convergence theorems for general step sizes require smoother unstable manifolds than can be proved using the cooperation hypothesis.

Benaïm and Hirsch then show that this result applies to any $m$-player generalized coordination game with stochastic fictitious play with optimal response determined as in the framework of [FK93], provided that the response map is smooth (which requires some noise). Generalized coordination games by definition have only two strategies per player, so the extension of these results to multi-strategy games was left open. At the time or writing, the final installment in the story of stochastic fictitious play is the extension by Hofbauer and Sandholm of the non-stochastic convergence results (for potential games, su-permodular games and games with an internal evolutionarily stable strategy) to the stochastic setting [HS02]. Forthcoming work of Benaïm, Hofbauer and Sorin [BHS05; BHS06] replaces the differential equation by a set valued **differential inclusion** in order to handle fictitious play with imperfect information or with discontinuous $F$.

## *4.6. Agent-based modeling*

In agent-based models, according to [Bon02], "A system is modeled as a collection of autonomous decision-making entities called agents, [with each] agent individually assessing its situation and making decisions on the basis of a set of rules." A typical example is a graph theoretic model, where the agents are vertices of a graph and at each time step, each agent chooses an action based on various characteristics of its neighbors in the graph; these actions, together with external sources of randomness, determine outcomes which may alter the characteristics of the agents. Stochastic replicator dynamics fall within this rubric, as do a number of the other processes already discussed. The boundaries are blurry, but this section is chiefly devoted to agent-based models from the social sciences, in which some sort of graph theoretic structure is imposed.

Analytic intractability is the rule rather than the exception for such models. The recent boom in agent-based modeling is probably due to the emergence of fast computers and of software platforms specialized to perform agent-based simulation. One scientific utility for such models is to give simple explanations for complex phenomena. Another motivation comes from psychology. Even in situations where people are capable of some kind of rational game-theoretic computation, evidence shows that actual decision mechanisms are often much more primitive. Brain architecture dictates that the different components of a decision are processed by different centers, with the responses then chemically or electrically superimposed (see for example [AHS05]). Three realistic components of decision making, captured better by agent-based models than by rational choice models are noted by Flache and Macy [FM02, page 633]:

- Players develop preferences for choices associated with better outcomes even though the association may be coincident, causally spurious, or superstitious.
- Decisions are driven by the two simultaneous and distinct mechanisms of reward and punishment, which are known to operate ubiquitously in humans.
- **Satisficing**, or persisting in a strategy that yields a positive but not optimal outcome, is common and indicates a mechanism of reinforcement rather than optimization.

Agent-based models now abound in a variety of social science disciplines, including psychology, sociology [BL03], public health [EL04], political science [OMH$^+$04]. The discussion here will concentrate on a few game-theoretic applications in which rigorous results have been obtained.

A number of recent analyses have centered on a two-player coordination game similar to Rousseau's **stag hunt**. Each player can choose to hunt rabbits or stags. The payoff is bigger for a stag but the stag hunt is successful only if both players hunt stag, whereas rabbit hunting is always successful. More generally, consider a payoff matrix as follows

$$\begin{pmatrix} (a,a) & (c,d) \\ (d,c) & (b,b) \end{pmatrix}. \tag{4.5}$$

When $a > d$ and $b > c$, the outcomes $(a,a)$ and $(b,b)$ are both Nash equilibria. Assume these inequalities, and without loss of generality, assume $a > b$. Then if $a - d > b - c$, the outcome $(b,b)$ is the risk-dominant equilibrium, whereas $(a,a)$ is always the unique Pareto-optimal equilibrium.

In 1993, Kandori, Mailath and Rob [KMR93] analyzed a very general class of evolutionary dynamics for populations of $N$ individuals associated with the two strategy types. The class included the following extreme version of stochastic replicator dynamics: each player independently with probability $1 - 2\epsilon$ changes type to whatever strategy type was most successful against the present population mix, and with probability $2\epsilon$ resets the type according to the result of independent fair coins. In the case of a game described by 4.5 they showed that the resulting Markov chain always converged to the risk-dominant equilibrium in the sense that the chain had a stationary measure $\mu_{N,\epsilon}$ satisfying:

**Theorem 4.4** ([KMR93, Theorem 3]). *As $\epsilon \to 0$ with $N$ fixed and sufficiently large, $\mu_{N,\epsilon}$ converges to the entire population playing the risk-dominant equilibrium.* □

PROOF: Assume without loss of generality that $a - d > b - c$, that is, that strategy 2 is risk-dominant. There is an embedded two-state Markov chain, where state 1 contains all populations where the proportion of type 1 players is at least $\alpha N$, and $\alpha(\epsilon)$ is the threshold for strategy 1 to be superior to strategy 2 against such a population. Due to $a - d > b - c$, we know $\alpha < N/2$. Going from state 2 to state 1 occurs exactly when there are at least $\alpha N$ "mutations" (types chosen by coin-flip) and going from state 1 to state 2 occurs when there are at least $\alpha N$ mutations. The ratio of the stationary measures of state 1 to state 2 goes to the ratio of these two probabilities, which goes to infinity. □

Unfortunately, the waiting time to get from either state to the other is exponential in $N \log(1/\epsilon)$, meaning that for many realistic parameter values, the population, if started at the sub-optimal equilibrium, does not have time to learn the better equilibrium. This many simultaneous mutations are as rare as all the oxygen molecules suddenly moving to the other side of the room (well not quite). Ellison [Ell93] proposes a variant. Let the agents be labeled by the integers modulo $N$, and for fixed $k < N/2$, let $i$ and $j$ be considered neighbors if their graph distance is at most $k$. Ellison's dynamics are the same as in [KMR93] except that each agent with probability $1 - 2\epsilon$ chooses the best play against the reference population consisting of that individual together with its $2k$ neighbors. The following result shows that when global interactions are replaced by local interactions, the population learns the optimal equilibrium much more rapidly.

**Theorem 4.5** ([Ell93, Theorem 3]). *For fixed $k$ and sufficiently small $\epsilon$, as $N \to \infty$ the expected time from any state to a state with most players of type 1 remains constant.*

PROOF: Let $j \leq k$ be such that $j$ out of $2k + 1$ neighbors of type 1 is sufficient to make strategy 1 optimal. Once there are $j$ consecutive players of type 1, the size of the interval of players of type 1 (allowing an $\epsilon$ fraction of errors) will tend to increase by roughly $2(k - j - \epsilon N)$ at each turn. The probability of the interval $r + 1, \ldots, r + j$ all turning to type 1 in one step is small but nonzero, so for sufficiently large $N$, such an interval arises immediately. □

The issue of how people might come to choose the superior $(a, a)$ in this case has been of longstanding concern to game theorists. In [SP00], a new evolutionary dynamic is introduced. A two-player game is fixed, along with a population of players labeled $1, \ldots, N$. Each player is initially assigned a strategy type. Positive weights $w(i, j, 1)$ are assigned as well, usually all equal to 1. The novel element to the model is the simultaneous evolution of network structure with strategy. Specifically, the network at time $t$ is given by the collection of weights $w(i, j, t)$ representing propensities for player $i$ to interact with player $j$ at time $t$. At each time step, each player $i$ chooses a partner $j$ independently at random with probabilities proportional to $w(i, j, t)$, then plays the game with the partner. After this, $w(i, j, t + 1)$ is set equal to $w(i, j, t) + u$ and $w(j, i, t + 1)$ is set equal to $w(j, i, t) + u'$, where $u$ and $u'$ are the respective utilities obtained by players $i$ and $j$. (Note that each player plays at least once in each round, but more than once if the player is chosen as partner by one of more of the other players.)

In their first model, Skyrms and Pemantle take the strategy type to be fixed and examine the results of evolving network structure.

**Theorem 4.6** ([SP00, Theorem 6]). *Consider a network of $2n$ players in Rousseau's Stag Hunting game given by the payoff matrix*

$$\begin{pmatrix} (1,1) & (0, 0.75) \\ (0.75, 0) & (0.75, 0.75) \end{pmatrix},$$

*with $2k > 0$ stag hunters and $2(n - k) > 0$ rabbit hunters. Under the above network evolution rules, with no evolution or mutation of strategies, as $t \to \infty$, the probability approaches 1 that all stag hunters choose stag hunters and all rabbit hunters choose rabbit hunters.*

PROOF: If $i$ is a stag hunter and $j$ is a rabbit hunter then $w(i, j, t)$ remains 1 for all time; hence stag hunters do not choose rabbit hunters in the limit. The situation is more complicated for $w(j, i, t)$, since rabbit hunters get reinforced no matter whom they choose or are chosen by. However, if $A$ denotes the set of stag hunters and $Z(j, t)$ denotes the probability $\sum_{i \in A} w(j, i, t) / \sum_i w(j, i, t)$ that $j$ will choose a stag hunter at time $t$, then it is not hard to find $\lambda, \mu > 0$ such that $\exp(\lambda Z(j, t) + \mu \log t)$ is a supermartingale, which implies that $Z(i, t) \to 0$ (in fact, exponentially fast in $\log t$). □

Further results via simulation show that when each agent after each round decides with a fixed probability $\epsilon > 0$ to switch to the strategy that is optimal against the present population, then all agents converge to a single type. However, it is random which type. When the evolution of strategy was slow (e.g., $\epsilon = 1/100$), the system usually found at the optimal equilibrium (everyone hunts stag) but when the evolution of strategy was more rapid (e.g., $\epsilon = 1/10$), the majority (78%) of the simulations resulted in the maximin equilibrium where everyone hunts rabbits. Evidently, more rapid evolution of strategy causes the system to mirror the stochastic replicator models in which the risk-dominant equilibrium is always chosen.

## 4.7. Miscellany

*Splines and interpolating curves*

Computer-aided drawing programs often provide interpolated curves. A finite sequence $x_0, \ldots, x_n$ of control points in $\mathbb{R}^d$ are specified, and a curve $\{f(t) : 0 \leq t \leq 1\}$ is generated which in some sense approximates the polygonal path $g(t)$ defined to equal $x_k + (nt - k)(x_{k+1} - x_k)$ for $k/n \leq t \leq (k+1)/n$. In many cases, the formula for producing $f$ is

$$f(t) = \sum_{k=0}^{n} B_{n,k}(t) x_k \,.$$

Depending on the choice of $\{B_{n,k}(t)\}$, one obtains some of the familiar blending curves: Bezier curves, B-splines, and so forth.

Goldman [Gol85] proposes a new family of blending functions. Consider a two-color Pólya urn with constant reinforcement $c \geq 0$, initially containing $t$ red balls and $1 - t$ black balls. Let $B_{n,k}(t)$ be the probability of obtaining exactly $k$ red balls in the first $n$ trials. The functions $\{B_{n,k}\}$ are shown to have almost all of the requisite properties for families of blending functions. In particular,

(i) $\{B_{n,k}(t) : k = 0, \ldots, n\}$ are nonnegative and sum to 1, implying that the interpolated curve is in the convex hull of the polygonal curve;
(ii) $B_{n,k}(t) = B_{n,n-k}(1 - t)$ implying symmetry under reversal;
(iii) $B_{n,k}(0) = \delta_{k,0}$ and $B_{n,k}(1) = \delta_{k,n}$, implying that the curve and polygon have the same endpoints (useful for piecing together curves);
(iv) $\sum_{k=0}^{n} k B_{n,k}(t) = nt$, implying that the curve is a line when $x_{k+1} - x_k$ is independent of $k$;
(v) The curve is less wiggly than the polygonal path: for any vector $\mathbf{v}$, the number of sign changes of $f(t) \cdot \mathbf{v}$ is at most the number of sign changes of $g(t) \cdot \mathbf{v}$
(vi) Given $P_0, \ldots, P_n$ there are $Q_0, \ldots, Q_{n+1}$ that reproduce the same curve $f(t)$ with the same parametrization;
(vii) Any segment $\{f(t) : a \leq t \leq b\}$ of the curve with control points $P_0, \ldots, P_n$ is reproducible as the entire curve corresponding to control points $Q_0, \ldots, Q_n$, where the parametrization may differ but $n$ remains the same.

There is of course an explicit formula for the polynomials $B_{n,k}(t)$. This generalizes the Bernstein polynomials, which are obtained when the reinforcement parameter $c$ is zero. However, the urn model pulls its weight in the sense that verification of many of the features is simplified by the urn model interpretation. For example, the first fact translates simply to the fact that for fixed $n$ and $t$ the quantities $B_{n,k}(t)$ are the probabilities of $k+1$ possible values of a random variable.

In a subsequent paper [Gol88a], Goldman goes on to represent the so-called Beta-spline functions of [Bar81] via a somewhat more complicated time-varying Friedman urn model. Classical B-splines have a similar representation, which has consequences for the closeness of approximations by B-splines and Bernstein polynomials [Gol88b].

*Image reconstruction*

An interesting application of a network of Pólya urns is described in [BBA99]. The object is to reconstruct an image, represented in a grid of pixels, each of which contains a single color from a finite color set $\{1, \ldots, k\}$. Some coherence of the image is presumed, indicating that pixels dissimilar to their neighbors are probably errors and should be changed to agree with their neighbors. Among the existing methods to do this are maximum likelihood estimators, Markov random field models with Gibbs-sampler updating, and smoothing via wavelets. Computation of the MLE may be difficult, the Gibbs sampler may converge too slowly, and wavelet computation may be time-consuming as well.

Banarjee *et al.* propose letting the image evolve stochastically via a network of urns. This is fast, parallelizable, and should capture the qualitative features of smoothing. The procedure is as follows. There is an urn for each pixel. Initially, urn $x$ contains $x(j)$ balls of color $j$, where

$$x(j) = \sum_{y \neq x} \frac{\delta(y, j)}{d(x, y)}$$

and $\delta(y, j)$ is one if pixel $y$ is colored $j$ and zero otherwise. In other words, the initial contents are determined by the empirical distribution of colors near $x$, weighted by inverse distance. Define a neighborhood structure: for each $x$ there is a set of pixels $N(x)$; this may for example be nearest neighbors or all pixels up to a certain distance from $x$. The update rule for urn $x$ is to sample from the combined urn of all elements of $N(x)$ and add a constant number $\Delta$ of balls of the sampled color to urn $x$. This may be done simultaneously for all $x$, sequentially, or by choosing $x$ uniformly at random. After a long time, the process halts and the output configuration is chosen by taking the plurality color at each pixel. The mathematical analysis is incomplete, but experimental data shows that this procedure outperforms a popular relaxation labeling algorithm (the urn scheme is faster and provides better noise reduction).

## 5. Reinforced random walk

In 1987 [Dia88] (see also [CD87]), Diaconis introduced the following process, known now as **edge-reinforced random walk** or **ERRW**. A walker traverses the edges of a finite graph. Initially any edge incident to the present location is equally likely, but as the process continues, the likelihood for the walker to choose an edge increases with each traversal of that edge, remaining proportional to the **weight** of the edge, which is one more than the number of times the edge has been traversed in either direction.

Formally, let $G := (V, E)$ be any finite graph and let $v \in V$ be the starting vertex. Define $X_0 = v$ and $W(e, 0) = 1$ for all $e \in E$. Inductively define $\mathcal{F}_n := \sigma(X_0, \ldots, X_n)$, $W(\{y, z\}, n) = W(\{y, z\}, n-1) + \mathbf{1}(\{X_{n-1}, X_n\} = \{y, z\})$, and let

$$\mathbb{P}(X_{n+1} = w | \mathcal{F}_n) = \frac{W(\{X_n, w\}, n)}{\sum_z W(\{X_n, z\}, n)}$$

if $w$ is a neighbor of $X_n$ and zero otherwise. The main result of [CD87] is that ERRW is a mixture of Markov chains, and that the edge occupation vector converges to a random limit whose density may be explicitly identified.

**Theorem 5.1** ([Dia88, (4.2)]). *Let $\{X_n\}$ be an ERRW on the finite graph $G = (V, E)$ beginning from the vertex $v_0$. Then $\{X_n\}$ is a mixture of Markov chains, meaning that there is a measure $\mu$ on transition probabilities $\{p(v, w) : \{v, w\} \in E\}$ such that*

$$\mathbb{P}(X_0 = v_0, \ldots, X_n = v_n) = \int p(v_0, v_1) \cdots p(v_{n-1}, v_n) \, d\mu.$$

*Furthermore, the weights $\mathbf{W} := \{W(e, n) : e \in E\}$ approach a random limit continuous with respect to Lebesgue measure on the simplex $\{\mathbf{W} : w(e) \geq 0, \sum_e w(e) = 1\}$ of sequences of nonnegative numbers indexed by $E$ and summing to 1. The density of the limit is given by*

$$C \prod_{e \in E} w(e)^{1/2} \prod_{v \in V} w(v)^{-(1+d(v))/2} w_{v_0}^{1/2} |A|^{1/2} \tag{5.1}$$

*where $w(v)$ denotes the sum of $w(e)$ over edges $e$ adjacent to $v$, $d(v)$ is the degree of $v$ and $A$ is the matrix indexed by cycles $C$ forming a basis for the homology group $H_1(G)$ with $A(C, C) := \sum_{e \in C} 1/w(e)$ and $A(C, D) = \sum_{e \in C \cap D} \pm 1/w(e)$ with a positive sign if $e$ has the same orientation in $C$ and $D$ and a negative sign otherwise.* □

This result is proved by invoking a notion of partial exchangeability [dF38], shown by [DF80] to imply that a process is a mixture of Markov chains[11]. The formula (5.1) is then proved by a direct computation. The computation was never written down and remained unavailable until a more general proof was

---

[11]When the process may not visit sites infinitely often, some care must be taken in deducing the representation as a mixture of Markov chains from partial exchangeability; see for example the recent paper of Merkl and Rolles [MR07].

published by Keane and Rolles [KR99]. The definition extends easily to ERRW on the infinite lattice $\mathbb{Z}^d$ and Diaconis posed the question of recurrence:

**Question 5.1.** *Does ERRW on $\mathbb{Z}^d$ return to the origin with probability 1?*

This question, still open, has provoked a substantial amount of study. Early results on ERRW and some of its generalizations are discussed in the next subsection; the following subsections concern two other variants: vertex-reinforced random walk and continuous time reinforced random walk on a graph. For further results on all sorts of ERRW models, the reader is referred to the short but friendly survey [MR06].

### 5.1. Edge-reinforced random walk on a tree

A preliminary observation is that ERRW on a directed graph may be represented by a network of Pólya urn processes. That is, suppose that $P(X_{n+1} = w \,|\, \mathcal{F}_n)$ is proportional to one plus the number of directed transits from $X_n$ to $w$. Then for each vertex $v$, the sequence of vertices visited after each visit to $v$ is distributed exactly as a Pólya urn process whose initial composition is one ball of color $w$ for each neighbor $w$ of $v$; as $v$ varies, these urns are independent. Formally, consider a collection of independent Pólya urns labeled by vertices $v \in V$, the contents of each of which are initially a single ball of color $w$ for each neighbor $w$ of $v$; let $\{X_{n,v} : n = 1, 2, \ldots\}$ denote the sequence of draws from urn $v$; then we may couple an ERRW $\{X_n\}$ to the independent urns so that $X_{n+1} = w \iff X_{s,v} = w$, where $s$ is the number of times $v$ has been visited at time $n$.

For the usual undirected ERRW, no such simple representation is possible because the probabilities of successive transitions out of $v$ are affected by which edges the path has taken coming into $v$. However, if $G$ is a tree, then the first visit to $v \neq v_0$ must be along the unique edge incident to $v$ leading toward $v_0$ and the $(n+1)^{st}$ visit to $v$ must be a reverse traversal of the edge by which the walk left $v$ for the $n^{th}$ time. This observation, which is the basis for Lemma 2.4, was used by [Pem88a] to represent ERRW on an infinite tree by an infinite collection of independent urns. In this analysis, the reinforcement was generalized from 1 to an arbitrary constant $c > 0$. The urn process corresponding to $v \neq v_0$ has initial composition $(1+c, 1, \ldots, 1)$ where the first component corresponds to the color of the parent of $v$, and reinforcement $2c$ each time. Recalling from (4.2) that such an urn is exchangeable with limit distribution that is Dirichlet with parameters $(1+c)/(2c), 1/(2c), \ldots, 1/(2c)$, one has a representation of ERRW on a tree by a mixture of Markov chains whose transition probabilities out of each vertex are given by picks from the specified independent Dirichlet distributions. This leads to the following phase transition result (see also extensions by Collevecchio [Col04; Col06a; Col06b]).

**Theorem 5.2** ([Pem88a, Theorem 1]). *There is a constant $c_0 \approx 4.29$ such that ERRW on an infinite binary tree is almost surely transient if $c < c_0$ and almost surely recurrent if $c > c_0$.* □

## *5.2. Other edge-reinforcement schemes*

The reinforcement scheme may be generalized in several ways. Suppose the transition probabilities out of $X_n$ at step $n$ are proportional not to the weights $w(\{X_n, w\}, n)$ incident to $X_n$ at time $n$ but instead to $F(w(\{X_n, w\}, n))$ where $F : \mathbb{Z}^+ \to \mathbb{R}^+$ is any nondecreasing function. Letting $a_n := F(n) - F(n-1)$, one might alternatively imagine that the reinforcement is $a_n$ on the $n^{th}$ time an edge is crossed (see the paragraph in Section 3.2 on ordinal dependence). Davis [Dav99] calls this a reinforced random walk of **sequence type**. A special case of this is when $a_1 = \delta$ and $a_n = 0$ for $n \geq 2$. This is called **once-reinforced random walk** for the obvious reason that the reinforcement occurs only once, and its invention is usually attributed to M. Keane. More generally, one might take the sequence to be different for every edge, that is, for each edge $e$ there is a nondecreasing function $F^e : \mathbb{Z}^+ \to \mathbb{R}^+$ and $\mathbb{P}(X_{n+1} = w \,|\, \mathcal{F}_n)$ is proportional to $F^e(w(e, n))$ with $e = \{X_n, w\}$.

It is easy to see that for random walks of sequence type on any graph, if $\sum_{n=1}^{\infty} 1/F(n) < \infty$ then with positive probability the sequence of choices out of a given edge will fixate. This extends to

**Theorem 5.3** ([Dav90, Theorem 3.2]). *Let $\{X_n\}$ be a random walk of sequence type on $\mathbb{Z}$. If $\sum_{n=1}^{\infty} 1/F(n) < \infty$ then $\{X_n\}$ almost surely is eventually trapped on a single edge. Conversely, if $\sum_{n=1}^{\infty} 1/F(n) = \infty$ then $\{X_n\}$ visits every vertex infinitely often almost surely.*

PROOF: Assume first that $\sum_{n=1}^{\infty} 1/F(n) < \infty$. To see that $\sup_n X_n < \infty$ with probability 1, it suffices to observe that for each $k$, conditional on ever reaching $k$, the probability that $\sup_n X_n = k$ is bounded below by $\prod_{n=1}^{\infty} F(n)/(1 + F(n))$ which is nonzero. The same holds for $\inf_n X_n$, implying finite range almost surely. To improve this to almost sure fixation on a single edge, Davis applies Herman Rubin's Theorem (Theorem 3.6) to show that the sequence of choices from each vertex eventually fixates. Conversely, if $\sum_{n=1}^{\infty} 1/F(n)$ is infinite, then each choice is made infinitely often from each vertex, immediately implying recurrence or converge to $\pm\infty$. The latter is ruled out by means of an argument based on the fact that the sum $M_n := \sum_{k=1}^{X_n} 1/F(w(\{j-1, j\}, n))$ of the inverse weights up to the present location is a supermartingale [Dav90, Lemma 3.0]. □

*Remark* 5.4. The most general ERRW considered in the literature appears in [Dav90]. There, the weights $\{w(e, n)\}$ are arbitrary random variables subject to $w(e, n) \in \mathcal{F}_n$ and $w(e, n+1) \geq w(e, n)$ with equality unless $e = \{X_n, \mathbf{X}_{n+1}\}$. The initial weights may be arbitrary as well, with the term **initially fair** used to denote all initial weights equal to 1. At this level of generality, there is no exchangeability and the chief techniques are based on martingales. Lemma 3.0 of [Dav90], used to rule out convergence to $\pm\infty$ is in fact proved in the context of such a general, initially fair ERRW.

When the graph is not a tree, many of the arguments become more difficult. Sellke [Sel94] extended the martingale technique to sequence-type ERRW on the $d$-dimensional integer lattice. Because of the bipartite nature of the graph, one must consider separately the sums $\sum_{n=1}^{\infty} 1/F(2n)$ and $\sum_{n=1}^{\infty} 1/F(2n+1)$. For convenience, let us assume these two both converge or both diverge.

**Theorem 5.5** ([Sel94, Theorems 1–3]). *If $\sum_{n=1}^{\infty} 1/F(n) < \infty$ then with probability one, the process is eventually trapped on a single edge. If $\sum_{n=1}^{\infty} 1/F(n) = \infty$, then with probability one, the range is infinite and each coordinate is zero infinitely often.* □

The proofs are idiosyncratic, based on martingales and Rubin's construction. It is noted that (*i*) the conclusion in the case $\sum_{n=1}^{\infty} 1/F(n) = \infty$ falls short of recurrence; and (*ii*) that the conclusion of almost sure trapping in the opposite case is specific to bipartite graphs, with the argument not generalizing to the triangular lattice, nor even to a single triangle! This situation was not remedied until Limic [Lim03] proved that for ERRW on a triangle, when $F(n) = n^\rho$ for $\rho > 1$, the walk is eventually trapped on a single edge. This was generalized in [LT06] to handle any $F$ with $\sum_{n=1}^{\infty} 1/F(n) < \infty$.

Because of the difficulty of proving results for sequence-type ERRW, it was thought that the special case of once-reinforced random walk might be a more tractable place to begin. Even here, no one has settled the question of recurrence versus transience for the two-dimensional integer lattice. The answer is known for a tree. In contrast to the phase transition in ordinary ERRW on a tree, a once-reinforced ERRW is transient for every $\delta > 0$ (in fact the same is true when "once" is replaced by "$k$ times"). This was proved for regular trees in [DKL02] and extended to Galton-Watson trees in [Die05].

The only other graph for which I am aware of an analysis of once-reinforced ERRW is the ladder. Let $G$ be the product of $\mathbb{Z}^1$ with $K_2$ (the unique connected two-vertex graph); the vertices are $\mathbb{Z} \times \{0,1\}$ and the edges connect neighbors of $\mathbb{Z}$ with the same $K_2$-coordinate or two vertices with the same $\mathbb{Z}$ coordinate. The following recurrence result was first proved by T. Sellke in 1993 in the more general context of allowing arbitrary vertical movement (cf. [MR05]).

**Theorem 5.6** ([Sel06, Theorem]). *For any $\delta > 0$, once-reinforced random walk on the ladder is recurrent.* □

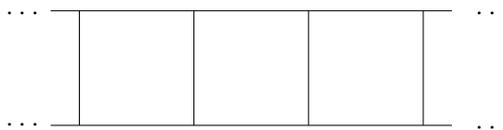

FIG 3. *The ladder graph of Theorem 5.6*

A recent result of Merkl and Rolles [MR05, Theorem 1.1] proves this for ERRW (as opposed to once-reinforced ERRW) as long as the ratio of the reinforcement parameter to the initial edge weights is less than $4/3$.

A slight variation on the definition of once-reinforced random walk gives the **excited random walk**, introduced in [BW03] and taken up in [Vol03; Zer06].

These results are too recent to have been included in this survey.

## 5.3. Vertex-reinforced random walk

Recall that the **vertex-reinforced random walk** (VRRW) is defined analogously to the ERRW except that in the equation (2.1) for choosing the next step, the edge occupation counts (2.2) are replaced by the vertex occupation counts (2.4).

This leads to entirely different behavior. Partial exchangeability is lost, so there is no representation as a random walk in a random environment. There are no obvious embedded urns. Moreover, an arbitrary occupation vector is unlikely to be evolutionarily stable. That is, suppose for some large $n$, the normalized occupation vector $\mathbf{X}_n$ whose components are the portion of the time spent at each vertex is equal to a vector $\mathbf{x}$. Let $\pi_\mathbf{x}$ denote the stationary measure for the Markov chain with transition probabilities $p(y, z) = \mathbf{x}_z / \sum_{z' \sim y} \mathbf{x}_{z'}$ which moves proportionally to the coordinate of $\mathbf{x}$ corresponding to the destination vertex. For $1 \ll k \ll n$, $\mathbf{X}_{n+k} = (1+o(1))\mathbf{X}_n$, so the proportion of the time in $[n, n+k]$ that the walk spends at vertex $y$ will be proportional to $\pi_\mathbf{x}(y)$. It is easy to see from this that $\{\mathbf{X}_n\}$ obeys a stochastic approximation equation (2.6) with $F(\mathbf{x}) = \pi_\mathbf{x} - \mathbf{x}$.

The analysis from here depends on the nature of the graph. The methods of Section 2.5 show that $X_n$ is an asymptotic pseudotrajectory for the flow $dX/dt = F(X)$, converging to an equilibrium point or orbit. There is always a Lyapunov function $V(\mathbf{x}) := \mathbf{x}^T A \mathbf{x}$ where $A$ is the incidence matrix of the underlying graph $G$. Therefore, equilibrium sets are sets of constancy for $V$ and any equilibrium point $p$ is a critical point for $V$ restricted to the face of the $(d-1)$-simplex containing $p$. Any attractor for the flow appears as a limit with positive probability, while linearly unstable orbits occur with probability zero. Several examples are given in [Pem88b].

**Example 5.1.** Let $G$ be the complete graph on $d$ vertices (with no self-edges). The zeros of $F$ are the centroids of faces. The global centroid $(1/d, \ldots, 1/d)$ is an attractor. Each other centroid is a permutation of some $(1/k, \ldots, 1/k, 0, \ldots, 0)$ and is easily seen to be linearly unstable [Pem88b, page 110]. It follows that $X_n \to (1/d, \ldots, 1/d)$ with probability 1.

**Example 5.2.** Let $G$ be a cycle of $d$ nodes for $d \geq 5$ (the smaller cases turn out to behave differently). The centroid $(1/d, \ldots, 1/d)$ is still an isolated equilibrium but for $d \geq 5$, it is linearly unstable. Although it was only guessed at the time this example appeared in [Pem88b], it follows from the nonconvergence theorems of [Pem90a; BH95] that the probability of convergence to the centroid is zero. The other equilibria are cyclic permutations of the points $(a, 1/2, 1/2 - a, 0, \ldots, 0)$ and certain convex combinations of these. It was conjectured in [Pem88b] and corroborated by simulation that the extreme points, namely cyclic permutations of $(a, 1/2, 1/2 - a, 0, \ldots, 0)$, were the only possible limits.

Taking $d \to \infty$ in the last example results in VRRW on the one-dimensional integer lattice. The analogous conjecture is that occupation measure for VRRW on $\mathbb{Z}$ converges to a translation of $\ldots, 0, 0, a, 1/2, 1/2 - a, 0, 0, \ldots$. It was shown [PV99, Theorem 1.1] that this happens with positive probability, leaving open the question of whether this occurs almost surely. Volkov strengthened and generalized this result. For infinite trees of bounded degree, VRRW gets trapped almost surely on a finite subtree [Vol01, Theorem 2]. In fact Volkov gives a graph theoretic definition of a **trapping subgraph** and shows that every finite graph has a trapping subgraph[12]. Volkov shows that on any locally finite connected graph without self-edges, if there is a trapping subgraph, $H$, then VRRW is trapped with positive probability on the union of $H$ and its neighbors. The neighbors of $H$ are visited with frequency zero, according to specified power laws [Vol01, Corollary 1]. Finally, Tarrès was able to close the gap from [PV99] for the one-dimensional lattice by proving almost sure trapping on an interval of exactly five vertices, with the conjectured power laws.

**Theorem 5.7** ([Tar04, Theorem 1.4]). *Let $\{X_n\}$ denote VRRW on $\mathbb{Z}$. With probability 1 there are (random) $k \in \mathbb{Z}, \alpha \in (0, 1)$ and $C_1, C_2 > 0$ such that the following occur.*

  (i) *The set of vertices visited infinitely often is $\{k-2, k-1, k, k+1, k+2\}$;*
 (ii) *The set of vertices visited with positive frequency is $\{k-1, k, k+1\}$ and these three frequencies have limits given respectively by $\alpha/2, 1/2$ and $(1-\alpha)/2$;*
(iii) *The occupation measure at $k-2$ is asymptotic to $C_1 n^\alpha$ and the occupation measure at $k+2$ is asymptotic to $C_2 n^{1-\alpha}$.*

□

### 5.4. An application and a continuous-time model

*Slime mold*

A mechanism by which simple organisms move in purposeful directions is called **taxis**. The organism requires a signal to govern such motion, which is usually something present in the environment such as sunlight, chemical gradient or particles of food.

Othmer and Stevens [OS97] consider instances in which the organism's response modifies the signal. In particular, Othmer and Stevens study myxobacteria: organisms which produce slime, over which it is then easier for bacteria to travel in the future. Aware of the work of Davis on ERRW [Dav90], they propose a stochastic cellular automaton to model the propagation of one or more bacteria. One of their goals is to determine what features of a model lead to stable aggregation of organisms; apparently previous such models have led to aggregates forming but then disbanding.

---

[12]In fact every graph of bounded degree has a trapping subgraph, though R. Thomas (personal communication) has found an infinite, locally finite graph with no trapping subgraph.

In the Othmer-Stevens model, the build-up in slime at the intersection points of the integer lattice is modeled by postulating that the likelihood for the organism to navigate to a given next vertex is one plus the number of previous visits to that site (by any organism). With one organism, this is just a VRRW, which they call a "simplified Davis' model", the simplification being to go from ERRW to VRRW. They allow a variable weight function $W(n) = \sum_{k=1}^{n} a_k$. On page 1047, Othmer and Stevens describe results from simulations of the "simplified" VRRW for a single particle. Their analysis of the simulations may be paraphrased as follows.

> If $F(n)$ grows exponentially, the particle ultimately oscillates between two points. If $F$ grows linearly with a small growth rate, the particle does not stay in a fixed finite region. These two results agree with the theoretical result, which is proven, however, only in one dimension. If the growth is linear with a large growth rate, results of the simulation are "no longer comparable to the theoretical prediction" but this is because the time for a particle to leave a fixed finite region increases with the growth rate of $F$.

Given what we know about VRRW, we can give a different interpretation of the simulation data. We know that VRRW, unlike ERRW, fixates on a finite set. The results of [Vol01] imply that for $\mathbb{Z}^2$ the fixation set has positive probability both of being a 4-cycle and of being a plus sign (a vertex and its four neighbors). All of this is independent of the linear growth rate. Therefore, the simulations with large growth rates do agree with theory: the particle is being trapped rather than exiting too slowly to observe. On the other hand, for small values of the linear reinforcement parameter, the particle must also be trapped in the end, and in this case it is the trapping that occurs too slowly to observe. The power laws in [Vol01, Corollary 1] and part *(iii)* of Theorem 5.7 give an indication of why the trapping may occur too slowly to observe.

Othmer and Stevens are ultimately concerned with the behavior of large collections of myxobacteria, performing a simultaneous VRRW (each particle at each time step chooses the next step independently, with probabilities proportional to the total reinforcement due to other any particle's visits to the destination site). They make the assumption that the system may be described by differential equations corresponding to the **mean-field limit** of the system, where the state is described by a density over $\mathbb{R}^2$. They then give a rigorous analysis of the mean field differential equations, presumably related to scaling limits of ERRW[13]. The mean-field assumption takes us out of the realm of rigorous mathematics, so we will leave Othmer and Stevens here, but in the end they are able to argue that stable aggregation may be brought about by the purely local mechanisms of reinforced random walk.

*A continuous-time reinforced jump process*

The next section treats a number of continuous-time models. I include the vertex-reinforced jump process in this section because it is a process on dis-crete space which does not involve a scaling limit and seems similar to the other models in this section.

---

[13] Analyses of these reaction-diffusion equations may be found in [LS97] as well.

The **vertex-reinforced jump process** (VRJP) is a continuous-time process on the one-dimensional lattice. From any site $x$, at time $t$, it jumps to each nearest neighbor $y$ at rate equal to one plus the amount of time, $L(y,t)$, that the process has spent at $y$. On a state space that keeps track of occupation measure as well as position, it is Markovian. The process is defined and constructed in [DV02] and attributed to W. Werner; because the jump rate at time $t$ is bounded by $2 + t$, the definition is completely routine. We may obtain a one-parameter family of reinforcement strengths by jumping at rate $C + L(y,t)$ instead of $1 + L(y,t)$.

The VRJP is a natural continuous-time analogue of VRRW. An alternative analogue would have been to keep the total jump rate out of $x$ at 1, the chance of a jump to $y = x \pm 1$ remaining proportional to the occupation measure at $y$. In fact the choice of variable jump rates decouples jumps to the left from jumps to the right, making the process more tractable. On any two consecutive sites $a$ and $a+1$, let $m(t)$ denote the occupation measure of $a+1$ the first time the occupation measure of $a$ is $t$. Then $m(t)/t$ is a martingale [DV02, Corollary 2.3], which implies convergence of the ratio of occupation measures at $a+1$ and $a$. Together with some computations, this leads to an exact characterization of the (unscaled) random limit normalized occupation measure.

**Theorem 5.8** ([DV02, Theorem 1.1]). *Let $L(n,t) := \int_0^t \mathbf{1}_{X_s = n}\, ds$ be the occupation measure at $n$ at time $t$. Then the limit $Y_n := \lim_{t \to \infty} t^{-1} L(n, t)$ exists for each $n$. Let $\{U_n\}$ be IID random variables with density*

$$\frac{\exp(-\frac{1}{2}(\sqrt{x} - \frac{1}{\sqrt{x}}))}{\sqrt{2\pi x^3}}.$$

*The collection $\{Y_n : n \in \mathbb{Z}\}$ is distributed as $\{W_n / \sum_{k=-\infty}^{\infty} W_k : n \in \mathbb{Z}\}$ where $W_0 = 1$, $W_n = \prod_{k=1}^n U_k$ if $n > 0$ and $W_k = \prod_{k=n}^{-1} U_k$ if $k < 0$.* □

This process may be defined on any locally finite graph. Limiting ratios of the occupation measure at neighboring vertices have the same description. On an infinite regular tree, this leads as in [Pem88a] to a transition between recurrence and transience, depending on the reinforcement parameter, $C$; see [DV04].

## 6. Continuous processes, limiting processes, and negative reinforcement

In this section we will consider continuous processes with reinforcement. Especially when these are diffusions, they might be termed "reinforced Brownian motion". Some of these arise as scaling limits of reinforced random walks, while others are defined directly. We then consider some random walks with negative reinforcement. The most extreme example is the self-avoiding random walk, which is barred from going where it has gone before. Limits of self-avoiding walks turn out to be particularly nice continuous processes.

### *6.1. Reinforced diffusions*

*Random walk perturbed at its extrema*

Recall the once-reinforced random walk of Section 5.2. This is a sequence-type ERRW with $F(0) = 1$ and $F(n) = 1 + \delta$ for any $n \geq 1$. The transition probabilities for this walk may be phrased as $\mathbb{P}(k, k+1) = 1/2$ unless $\{X_n\}$ is at its maximum or minimum value, in which case $\mathbb{P}(k, k+1) = 1/(2+\delta)$ or $(1+\delta)/(2+\delta)$ respectively.

If such a process has a scaling limit, the limiting process would evolve as a Brownian motion away from its left-to-right maxima and minima, plus some kind of drift inwards when it is at a left-to-right extremum. This inward drift might come from a local time process, but constructions depending on local time processes involve considerable technical difficulty (see, e.g., [TW98]). An alternate approach is an implicit definition that makes use of the maximum or minimum process, recalling the way a reflecting Brownian motion may be constructed as the difference, $\{B_t - B_t^\#\}$, between a Brownian motion and its minimum process.

Let $\alpha, \beta \in (-\infty, 1)$ be fixed, let $g^*(t) := \sup_{0 \leq s \leq t} g(s)$ denote the maximum process of a function $g$ and let $g^\#(t) := \inf_{0 \leq s \leq t} g(s)$ denote its minimum process. Carmona, Petit and Yor [CPY98] examine the equation

$$g(t) = f(t) + \alpha g^*(t) + \beta g^\#(t), \qquad\qquad t \geq 0. \qquad (6.1)$$

They show that if $f$ is any continuous function vanishing at 0, then there is a unique solution $g(t)$ to (6.1), provided that $\rho := |\alpha\beta/((1-\alpha)(1-\beta))| < 1$. If $f$ is the sample path of a Brownian motion, then results of [CPY98] imply that the solution $Y_t := g(t)$ to (6.1) is adapted to the Brownian filtration. It is a logical candidate for a "Brownian motion perturbed at its extrema".

In 1996, Burgess Davis [Dav96] showed that the Carmona-Petit-Yor process is in fact the scaling limit of the once-reinforced random walk. His argument is based on the property that the map taking $f$ to $g$ in (6.1) is bounded: $|g_1 - g_2| \leq C|f_1 - f_2|$. The precise statement is as follows. Let $\alpha = \beta = -\delta$. The process $g$ will be well defined, since $\rho = |\delta/(1-\delta)|^2 < 1$.

**Theorem 6.1** ([Dav96, Theorem 1.2]). *Let $\{X_n : n \geq 0\}$ be the once-reinforced random walk with parameter $\delta > 0$. Let $\{Y_t : t \geq 0\}$ solve (6.1) with $\alpha = \beta = -\delta$ and $f$ a Brownian motion. Then*

$$\{n^{-1/2} X_{\lfloor nt \rfloor} : t \geq 0\} \xrightarrow{\mathcal{D}} \{Y_t : t \geq 0\}$$

*as $n \to \infty$.* □

*Drift as a function of occupation measure*

Suppose one wishes to formulate a diffusion that behaves like a Brownian motion, pushed according to a drift that depends in some natural way on the past, say through the occupation measure. There are a multitude of ways to do this. One way, suggested by Durrett and Rogers [DR92], is to choose a function $f$ and let the drift of the diffusion $\{X_t\}$ be given by $\int_0^t f(X_t - X_s)\, ds$. If $f$ is Lipschitz then there is no trouble in showing that the equation

$$X_t = B_t + \int_0^t ds \int_0^s du\, f(X_s - X_u)$$

has a pathwise unique strong solution.

Durrett and Rogers were the first to prove anything about such a process, but they could not prove much. When $f$ has compact support, they proved that in any dimension there is a nonrandom bound $\limsup_{t \to \infty} |X_t|/t \leq C$ almost surely, and that in one dimension when $f \geq 0$ and $f(0) > 0$, then $X_t/t \to \mu$ almost surely for some nonrandom $\mu$. The condition $f \geq 0$ was weakened by [CM96] to be required only in a neighborhood of 0. Among Durrett and Rogers' conjectures are that if $f$ is a compactly supported odd function with $xf(x) \geq 0$, then $X_t/t \to 0$ almost surely. Their reasoning is that the process should behave roughly like a negatively once-reinforced random walk. It sees only the occupation in an interval, say $[X_t - 1, X_t + 1]$, drifting linearly to the right for a while due to the imbalance in the occupation, until diffusive fluctuations cause it to go to the left of its maximum. Now it should get pushed to the left at roughly linear rate until it suffers another reversal. They were not able to make this rigorous. However, taking the support to zero while maintaining $\int_{-\infty}^0 f(x)\, dx = c$ gives a very interesting process about which Toth and Werner were able to obtain results (see Section 6.3 below).

Cranston and Le Jan [CLJ95] take up this model in two special cases. When $f(x) = -ax$ with $a > 0$, there is a restoring force equal to the total moment of the occupation measure about the present location. The restoring force increases without bound, so it may not be too surprising that $X_t$ converges almost surely to the mean of the limiting occupation measure.

**Theorem 6.2** ([CLJ95, Theorem 1]). *Let $a > 0$ and set $f(x) = ax$. Then there is a random variable $X_\infty$ such that $X_t \to X_\infty$ almost surely and in $L^2$.*

PROOF: This may be derived from the fact that the stochastic differential equation

$$dX_t = dB_t - \left[a \int_0^t (X_t - X_u)\, du\right] dt$$

has the (unique) strong solution

$$X_t = \int_0^t h(t,s) dB_t$$

with $h(t,s) = 1 - ase^{as^2/2} \int_s^t e^{-au^2/2}\, du$. □

The other case they consider is $f(x) = -a\,\mathrm{sgn}(x)$. It is not hard to show existence and uniqueness despite the discontinuity. This time, the restoring force is toward the median rather than the mean, but otherwise the same result, $X_t \to X$ should and does hold [CLJ95, Theorem 2]. This was extended to higher dimensions by the following theorem of Raimond.

**Theorem 6.3** ([Rai97, Theorem 1]). *Let $B_t$ be $d$-dimensional Brownian motion and fix $\sigma > 0$. Then the process $\{X_t\}$ that solves $X(0) = 0$ and*

$$dX_t = dB_t - \sigma \int_0^t \frac{X_t - X_s}{||X_t - X_s||} \, ds \, dt$$

*converges almost surely.* □

*Drift as a function of normalized occupation measure*

The above diffusions have drift terms that are additive functionals of the full occupation measure $\mu_t := \mu_0 + \int_0^t \delta_{X_s} \, ds$. The papers that analyze this kind of diffusion are [DR92; CLJ95; Rai97; CM96]; see also [NRW87]. In a series of papers [BLR02; BR02; BR03; BR05], Benaïm and Raimond (and sometimes Ledoux), consider diffusions whose drift is a function of the normalized occupation measure $\pi_t := t^{-1}\mu_t$. Arguably, this is closer in spirit to the reinforced random walk. Another difference in the direction taken by Benaïm and Raimond is that their state space is a compact manifold without boundary. This sets it apart from continuum limits of reinforced random walks on $\mathbb{Z}^d$ (not compact) or limits of urn processes on the $(d-1)$-simplex (has boundary).

The theory is a vast extension of the dynamical system framework discussed in Section 2.5. To define the object of study, let $M$ be a compact Riemannian manifold. There is a Riemannian probability measure, which we call simply $dx$, and a standard Brownian motion defined on $M$ which we call $B_t$. Let $V : M \times M$ be a smooth "potential" function and define the function $V\mu$ by

$$V\mu(y) = \int V(x, y) \, d\mu(x) \,.$$

The additive functional of normalized occupation measure is always taken to be $V\pi_t = t^{-1}V\mu_t$; thus the drift at time $t$ should be $\nabla(V\pi_t)$. Since $V\mu(\cdot) = \int V(x, \cdot) \, d\mu(x)$, we may write the stochastic differential equation as in [BLR02]:

$$dX_t = dB_t - \frac{1}{t}\left[\int_0^t \nabla V(X_s, X_t) \, ds\right] dt \,. \tag{6.2}$$

A preliminary step establishes the existence of the process $\{X_t\}$ from any starting point, and including the possibility of any starting occupation measure [BLR02, Proposition 2.5]. Simultaneously, this defines the occupation measure process $\{\mu_t\}$ and normalized occupation measure process $\{\pi_t\}$. When $t$ is large, $\pi_{t+s}$ will remain near $\pi_t$ for a while. As in the dynamical system and stochastic approximation framework, the next step is to investigate what happens if one fixes the drift for times $t + s$ at $-\nabla(V\pi_t)$. A diffusion on $M$ with drift $-\nabla f$ has an invariant measure whose density may be described explicitly as

$$\frac{\exp(-f(x))}{\int_M \exp(-f(z))\, dz}\,.$$

This leads us to define a function $\Pi$ that associates with each measure $\mu$ the density of the stationary measure for Brownian motion with potential function $V\mu$:

$$\Pi(\mu) := \frac{\exp(-V\mu(x))}{\int_M \exp(-V\mu(z))\, dz}\, dx\,. \qquad (6.3)$$

The process $\{\pi_t\}$ evolves stochastically. Taking a cue from the framework of Section 2.5, we compute the deterministic equation of mean flow. If $t \gg \delta t \gg 1$, then $\pi_{t+\delta t}$ should be approximately $\pi_t + \frac{\delta t}{t}(\Pi(\pi_t) - \pi_t)$. Thus we are led to define a vector field on the space of measures on $M$ by

$$F(\mu) := \Pi(\mu) - \mu\,.$$

A second preliminary step, carried out in [BLR02, Lemma 3.1], is that this vector field is smooth and induces a flow $\Phi_t$ on the space $\mathcal{P}(M)$ of probability measures on $M$ satisfying

$$\frac{d}{dt}\Phi_t(\mu) = F(\Phi_t(\mu))\,.$$

As with stochastic approximation processes, one expects the trajectories of the stochastic process $\pi_t$ to approximate trajectories of $\Phi_t$. One expects convergence of $\pi_t$ to fixed points or closed orbits of the flow, positive probability of convergence to isolated sinks, and zero probability of convergence to unstable equilibria. A good part of the work accomplished in the sequence of papers [BLR02; BR02; BR03; BR05] is to extend results on asymptotic pseudotrajectories in $\mathbb{R}^d$ to prove these convergence and nonconvergence results in the space of measures on $M$. One link in the chain that does not need to be extended is that Theorem 2.14 (asymptotic pseudotrajectories have chain transitive limits), which is already valid in a general metric space. Benaïm *et al.* then go on to prove the following results. The proof of the first is quite technical and occupies Section 5 of [BLR02].

**Theorem 6.4** ([BLR02, Theorem 3.6])**.** *The flow $\pi_{e^t}$ is an asymptotic pseudotrajectory for the flow $\Phi_t$.* □

**Corollary 6.5.** *The limit set of $\{\pi_t\}$ is almost surely an invariant chain-recurrent set containing no proper attractor.*

**Theorem 6.6** ([BR05, Theorem 2.4])**.** *Suppose that $V$ is symmetric, that is, $V(x, y) = V(y, x)$. With probability 1, the limit set of the process $\{\pi_t\}$ is a compact connected subset of the fixed points of $\Pi$ (that is, the zero set of $F$).*

PROOF: Define the **free energy** of a strictly positive $f \in L^2(dx)$ by

$$J(f) := J_V(f) := \frac{1}{2}\langle Vf, f\rangle + \langle f, \log f\rangle$$

where $Vf$ denotes the potential $Vf(y) = \int V(x,y)f(x)\,dx$ and $\langle f, g \rangle$ denotes $\int f(x)g(x)\,dx$. Next, verify that $J$ is a Lyapunov function for the flow $\Phi_t$ ([BR05, Proposition 4.1]) and that $F(\mu) = 0$ if and only if $\mu$ has a density $f$ and $f$ is a critical point for the free energy, i.e., $\nabla J(F) = 0$ ([BR05, Proposition 2.9]). The result then follows with a little work from Theorem 6.4 and the general result that an asymptotic pseudotrajectory is a chain transitive set. □

**Corollary 6.7** ([BR05, Corollary 2.5]). *If, in addition, the zero set of $F$ contains only isolated points then $\pi_t$ converges almost surely.* □

The next two results are proved in a similar manner to the proofs of the convergence and nonconvergence results Theorem 2.16 and Theorem 2.9, though some additional infrastructure must be built in the infinite-dimensional case.

**Theorem 6.8** ([BR05, Theorem 2.24]). *If $\pi^*$ is a sink then $\mathbb{P}(\pi_t \to \pi^*) > 0$.* □

The nonconvergence results, as well as criteria for existence of a sink, the following definition is very useful.

**Definition 6.9** (Mercer kernel). *Say that $V$ is a Mercer kernel if $\langle Vf, f \rangle \geq 0$ for all $f \in L^2(dx)$.*

While the assumption of a symmetric Mercer kernel may appear restrictive, it is shown [BR05, Examples 2.14–2.20] that many classes of kernels satisfy this, including the transition kernel for any reversible Markov semi-group, any even function of $x - y$ on the torus $T^n$ that has nonnegative Fourier coefficients, any completely monotonic function of $||x - y||^2$ for a manifold embedded in $\mathbb{R}^n$. and any $V$ represented as $V(x,y) = \int_E G(\alpha, x)G(\alpha, y)\,d\nu(\alpha)$ for some space $E$ and measure $\nu$ (this last class is in fact dense in the set of Mercer measures). The most important fact about Mercer kernels is that they are strictly convex.

**Lemma 6.10** ([BR05, Theorem 2.13]). *If $V$ is Mercer then $J$ is strictly convex, hence has a unique critical point $f$ which is a global minimum.*

**Corollary 6.11.** *If $V$ is Mercer then the process $\pi_t$ converges almost surely to the measure $f\,dx$ where $f$ minimizes the free energy, $J$.* □

PROOF OF LEMMA: The second derivative $D^2 J$ is easily computed [BR05, Proposition 2.9] to be

$$D^2 f(u, v) \langle Vu, v \rangle_{dx} + \langle u, v \rangle_{(1/f)\,dx}.$$

The second term is always positive definite, while the first is nonnegative definite by hypothesis. □

The only nonconvergence result they prove requires a hypothesis involving Mercer kernels.

**Theorem 6.12** ([BR05, Theorem 2.26]). *If $\pi^*$ is a quadratically nondegenerate zero of $F$ with at least one positive eigenvalue, and if $V$ is the difference of Mercer kernels, then $\mathbb{P}(\pi_t \to \pi^*) = 0$.* □

A number of examples are given but perhaps the most interesting is one where $V$ is not symmetric. It is possible that there is no Lyapunov function and that the limit set of $\pi_t$, which must be an asymptotic pseudotrajectory, may be a nontrivial orbit. In this case, one expects that $\mu_t$ should precess along the orbit at logarithmic speed, due to the factor of $1/t$ in the mean differential equation $d\pi_t/dt = (1/t)F(\pi_t)$.

**Example 6.1.** Let $M = S^1$ and for fixed $c$ and $\phi$ let

$$V(\theta_1, \theta_2) = 2c \cdot \cos(\theta_1 - \theta_2 + \phi).$$

When $\phi$ is not 0 or $\pi$ this is not symmetric. A detailed trigonometric analysis shows that when $c \cdot \cos \phi \geq -1/2$, the unique invariant set for $\Phi_t$ is Lebesgue measure, $dx$, and hence that $\pi_t \to dx$ almost surely.

Suppose now that $c \cdot \cos \phi < -1/2$. If $\phi = 0$ then the critical points of the free energy function are a one-parameter family of zeros of $F$ with densities $g_\theta := c_1(c) e^{c_2(c) \cdot \cos(x - \theta)}$. It is shown in [BLR02, Theorem 1.1] that $\pi_t \to g_Z$ almost surely, where $Z$ is a random variable. The same holds when $\phi = \pi$.

When $\phi \neq 0, \pi$ then things are the most interesting. The forward limit set for $\{\pi_t\}$ under $\Phi$ consists of the unstable equilibrium point $dx$ (Lebesgue measure) together with a periodic orbit $\{\rho_\theta : \theta \in S^1\}$, obtained by averaging $g_\theta$ while moving with logarithmic speed. To rule out the point $dx$ as a limit for the stochastic process 6.2 would appear to require generalizing the nonconvergence result Theorem 2.17 to the infinite-dimensional setting. It turns out, however, that the finite-dimensional projection $\mu \mapsto \int_{S^1} x \, d\mu$ maps the process to a stochastic approximation process in the unit disk, that is, the evolution of $\int_{S^1} x \, d\mu$ depends on $\mu$ only through $\int_{S^1} x \, d\mu$. For the projected process, 0 is an unstable equilibrium, whence $dx$ is almost surely not a limit point of $\{\pi_t\}$. By Corollary 6.5, the limit set of the process is the periodic orbit. In fact there is a random variable $Z \in S^1$ such that

$$||\pi_t - \rho_{\log t + Z}|| \to 0.$$

This precise result relies on shadowing theorems such as [Ben99, Theorem 8.9].

### 6.2. Self-avoiding walks

A path of finite length on the integer lattice is said to be **self-avoiding** if its vertices are distinct. Such paths have been studied in the context of polymer chemistry beginning with [Flo49], where nonrigorous arguments were given to show that the diameter of a polymer chain of length $n$ in three-space would be of order $n^\nu$ for some $\nu$ greater than the value of $1/2$ predicted by a simple random walk model. Let $\Omega_n$ denote the set of self-avoiding paths in $\mathbb{Z}^d$ of length $n$ starting from the origin. Surprisingly, good estimates on the number of such paths are still not known. Hammersley and Morton [HM54] observed that $|\Omega_n|$ is sub-multiplicative: concatenation is a bijection between $\Omega_j \times \Omega_k$ and a set containing $\Omega_{j+k}$. It follows that $|\Omega_n|^{1/n}$ converges to $\inf_k |\Omega_k|^{1/k}$. The **connective constant** $\mu = \mu_d$, defined to be the value of this limit in $\mathbb{Z}^d$, is not known, though rigorous estimates place $\mu_2 \in [2.62, 2.70]$ and nonrigorous estimates claim great precision. It is not known though widely believed that in any dimension, $|\Omega_{n+1}|/|\Omega_n| \to \mu_d$; for closed loops Kesten [Kes63] did show that $|\Omega_{n+2}|/|\Omega_n| \to \mu^2$.

Let $U_n$ denote the uniform measure on $\Omega_n$. Given that the cardinality of $\Omega_n$ is poorly understood, it is not surprising that $U_n$ is also poorly understood. In dimensions five and higher, a substantial body of work by Hara and Slade has established the convergence under rescaling of $U_n$ to Brownian motion, convergence of $|\Omega_{n+1}|/|\Omega_n|$, and values of several exponents and constants. Their technique is to use asymptotic expansions known as **lace expansions**, based on numbers of various sub-configurations in the path. See [MS93] for a comprehensive account of work up to 1993 or Slade's piece in the *Mathematical Intelligencer* [Sla94] for a nontechnical overview.

In dimensions 2, 3 and 4, very little is rigorously known. Nevertheless, there are many conjectures, such as the existence and supposed values of diffusion exponent $\nu = \nu_d$ for which the $U_n$-expected square distance between the endpoints of the path (usually denoted $\langle R_n^2 \rangle$) is of order $n^{2\nu}$. Absent rigorous results, the measure $U_n$ has been investigated by simulation, but even that is difficult. The exponential growth of $U_n$ prevents sampling for $U_n$ in any direct way once $n$ is of order, say, 100.

Various Monte Carlo sampling schemes have been proposed. Beretti and Sokal [BS85] suggest a Markov Chain Monte Carlo algorithm, each step of which either extends the or retracts the path by one edge. Adjusting the relative probabilities of extension and retraction produces a Markov chain whose stationary distribution approximates a mixture of the measures $U_n$ and which approaches this distribution in polynomial time, provided certain conjectures hold and parameters have been correctly adjusted. Randall and Sinclair take this a step further, building into the algorithm foolproof tests of both of these provisions [RS00].

Of relevance to this survey are dynamic reinforcement schemes to produce self-avoiding or nearly self-avoiding random walks. It should be mentioned that there is no consensus on what measure should properly be termed the infinite self-avoiding random walk. If $U_n$ converges weakly then the limit is a candidate for such a walk. Two other ideas, discussed below, are to make the self-avoiding constraint soft, then take a limit, and to get rid of self-intersection by erasing loops as they form.

*'True' self-avoiding random walk*

For physicists, it is natural to consider the constraint of self-avoidance to be the limit of imposing a finite penalty for each self-intersection. In such a formulation, the probability of a path $\gamma$ is proportional to $e^{-\beta \cdot H(\gamma)}$ where the energy $H(\gamma)$ is the sum of the penalties.

A variant on this is to develop the walk dynamically via

$$\mathbb{P}(X_{n+1} = y \,|\, X_n = x, \mathcal{F}_n) = \frac{e^{-\beta \cdot N(y,n)}}{\sum_{z \sim x} e^{-\beta \cdot N(z,n)}}$$

where $N(z, n)$ is the number of visits to $z$ up to time $n$. This does not yield the same measure as the soft-constraint ensemble, but it has the advantage that it extends to a measure on infinite paths. Such a random walk was first considered by [APP83] and given the unfortunate name **true self-avoiding walk**. For finite inverse temperature $\beta$, this object is nontrivial in one dimension as well as in higher dimensions, and most of what is rigorously known pertains to one dimension. Tóth [Tót95] proves a number of results. His penalty function counts the number of pairs of transitions across the same edge, rather than the number of pairs of times the walk is at the same vertex, but is otherwise the same as in [APP83].

In the terminology of this survey, we have an ERRW [Tót95] or VRRW [APP83] of sequence type, with sequence $F(n) = e^{-\beta n}$. Tóth calls this exponential self-repulsion. In a subsequent paper [Tót94], the dynamics are generalized to subexponential self-repulsion $F(n) = e^{-\beta n^\kappa}$, with $0 < \kappa < 1$. The papers [Tót96; Tót97] then consider polynomial reinforcement $F(n) = n^\alpha$. When $\alpha < 0$ this is self-repulsion and when $\alpha > 0$ it is self-attraction. The following results give a glimpse into this substantial body of work, concentrating on the case of self-repulsion. An overview, which includes the case of self-attraction, may be found in the survey [Tót99]. For technical reasons, instead of $X_N$ in the first result, a random stopping time $\theta(N) = \theta(\lambda, N)$ is required. Define $\theta(N)$ to be a geometric random variable with mean $\lambda N$, independent of all other variables.

**Theorem 6.13** (see [Tót99, Theorem 1.4]). *Let $\{X_n\}$ be a sequence-type ERRW with sequence $F(n)$ equal to one of the following, and define the constant $\nu$ in each case as shown.*

1. $F(n) = \exp(-\beta n)$ , $\nu = \frac{2}{3}$

2. $F(n) = \exp(-\beta n^\kappa)$ , $\nu = \frac{\kappa+1}{\kappa+2}$

3. $F(n) = n^{-\alpha}$ , $\nu = \frac{1}{2}$

*where $\alpha$ is a positive constant and $0 < \kappa < 1$. There is a one-parameter family of distribution functions $\{G_\lambda(t) : \lambda > 0\}$ such that as $N \to \infty$,*

$$N^{-\nu} \mathbb{P}(X_{\theta(N)} < x) \to G_\lambda(x) \,.$$

*The limit distribution $G_\lambda$ is not Gaussian even when $\nu = 1/2$.* □

Tóth also proves a Ray-Knight theorem for the local time spent on each edge. As expected the time scaling is $N^{-\gamma}$ where $\gamma = (1-\nu)/\nu$.

*Loop-erased random walk*

Lawler [Law80] introduced a new way to generate a random self-avoiding path. Assume the dimension $d$ is at least 3. Inductively, we suppose that at time $n$, a self-avoiding walk from the origin $\gamma_n := (x_0, x_1, \ldots, x_k) \in \Omega_k$ has been chosen. Let $X_{n+1}$ be chosen uniformly from the neighbors of $x_k$, independently of what has come before. At time $n+1$, if $X_{n+1}$ is distinct from all $x_j$, $0 \leq j \leq k$, then $\gamma_{n+1}$ is taken to be $(x_0, \ldots, x_k, X_{n+1})$. If not, then $\gamma_{n+1}$ is taken to be $(x_0, \ldots, x_r)$ for the unique $r \leq n$ such that $x_r = X_{n+1}$ (we allow $\gamma_{n+1}$ to become the empty sequence if $r = 0$). In other words, the final loop in the path $(x_0, \ldots, x_r, x_{r+1}, \ldots, X_{n+1})$ is erased. In dimension three and higher, $|X_n| \to \infty$, and hence for each $k$ the first $k$ steps of $\gamma_n$ are eventually constant. The limiting path $\gamma$ is therefore well defined on a set of probability 1 and is a deterministic function, the **loop-erasure** of the simple random walk path $X_0, X_1, X_2, \ldots$, denoted $\mathsf{LE}(\mathbf{X})$. The loop-erased random walk measure, $\mathsf{LERW}$ is defined to be the law of $\gamma$.

The process $\gamma = \mathsf{LE}(\mathbf{X})$ seems to have little to do with reinforcement until one sees the following alternate description. Let $\{Y_n : n \geq 0\}$ be defined inductively by $Y_0 = 0$ and

$$\mathbb{P}(Y_{n+1} = y \mid Y_0, \ldots, Y_n) = \frac{h(y)}{\sum_{z \sim Y_n} h(z)} \qquad (6.4)$$

where $h(z)$ is the probability that a simple random walk beginning at $z$ avoids $\{Y_0, \ldots, Y_n\}$ forever, with $h(z) := 0$ if $z = Y_k$ for some $k \leq n$. Lawler observed [Law91, Proposition 7.3.1] that $\{Y_n\}$ has law $\mathsf{LERW}$. Thus one might consider $\mathsf{LERW}$ to be an infinitely negatively reinforced VRRW that sees the future. Moreover, altering (6.4) by conditioning on avoiding the past for time $M$ instead of forever, then letting $M \to \infty$ gives a definition of $\mathsf{LERW}$ in two dimensions that agrees with the loop-erasing construction when both are stopped at random times.

The law of $\gamma = \mathsf{LE}(\mathbf{X})$ is completely different from the laws $U_n$ and their putative limits, yet has some very nice features that make it worthy of study. It is time reversible, so for example the loop erasure of a random walk from $a$ conditioned to hit $b$ and stopped when it hits $b$ has the same law if $a$ and $b$ are switched. The loop-erased random walk on an arbitrary graph is also intimately related to an algorithm of Aldous and Broder [Ald90] for choosing a spanning tree uniformly. In dimensions five and above, $\mathsf{LERW}$ behaves the same way as the self-avoiding measure $U_n$, rescaling to a Brownian motion, but in dimensions 2, 3 and 4, it has different connectivity and diffusion exponents from $U_n$.

### 6.3. Continuous time limits of self-avoiding walks

Both the 'true' self-avoiding random walk and the loop-erased random walk have continuous limiting processes that are very pretty. The chance to spend a few paragraphs on each of these was a large part of my reason for including the entire section on negative reinforcement.

*The 'true' self-repelling motion*

The true self-avoiding random walk with exponential self-repulsion was shown in Theorem 6.13 (part 1) to have a limit law for its time-$t$ marginal. In fact it has a limit as a process. Most of this is shown in the paper [TW98], with a key tightness result added in [NR06]. Some properties of this limit process $\{X_t\}$ are summarized as follows. In particular, having 3/2-variation it is not a diffusion.

- The process $\{X_t\}$ has continuous paths.
- It is recurrent.
- It is self-similar:
$$\{X_t\} \stackrel{\mathcal{D}}{=} \{\alpha^{-2/3} X_{\alpha t}\}.$$
- It has non-trivial local variation of order 3/2.
- The occupation measure at time $t$ has a density; this may be called the local time $L_t(x)$.
- The pair $(X_t, L_t(\cdot))$ is a Markov process.

To construct this process and show it is the limit of the exponentially repulsive true self-avoiding walk, Tóth and Werner rely on the Ray-Knight theory developed in [Tót95]. While technical statements would involve too much notation, the gist is that the local time at the edge $\{k, k+1\}$ converges under re-scaling, not only for fixed $k$ but as a process in $k$. A strange but convenient choice is to stop the process when the occupation time on an edge $z$ reaches $m$. The joint occupations of the other edges $\{j, j+1\}$ then converge, under suitable rescaling, to a Brownian motion started at time $z$ and position $m$ and absorbed at zero once the time parameter is positive; if $z < 0$ it is reflected at zero until then. When reading the previous sentence, be careful, as Ray-Knight theory has a habit of switching space and time.

Because this holds separately for each pair $(z, m) \in \mathbb{R} \times \mathbb{R}^+$, the limiting process $\{X_t\}$ may be constructed in the strong sense by means of coupled coalescing Brownian motions $\{B_{z,m}(t) : t \geq z\}_{z \in \mathbb{R}, m \in \mathbb{R}^+}$. These coupled Brownian motions are jointly limits of coupled simple random walks. On this level, the description is somewhat less technical, as follows.

Let $V_e$ denote the even vertices of $\mathbb{Z}^2 \times \mathbb{Z}^+$. For each $(z, m) \in V_e$, flip an independent fair coin to determine a single directed edge from $(z, m)$ to $(z+1, m\pm 1)$; the exception is when $m = 1$; then for $z < 0$ there is an edge $\{(z, 1), (z+1, 2)\}$ while for $z \geq 0$ there is a v-shaped edge $\{(z, 1), (z+1, 0), (z+2, 1)\}$. Traveling rightward, one sees coalescing simple random walks, with absorption at zero once time is positive. A picture of this is shown. If one uses the even sites and travels leftward, one obtains a dual, distributed as a reflection (in time) of the original coalescing random walks. The complement of the union of the coalescing random walks and the dual walks is topologically a single path. Draw a polygonal path down the center of this path: the $z$-values when the center line crosses an integer level form a discrete process $\{Y_n\}$.

This process $\{Y_n\}$ is a different process from the true self-avoiding walk we started with, but it has some other nice descriptions, discussed in [TW98, Sec-

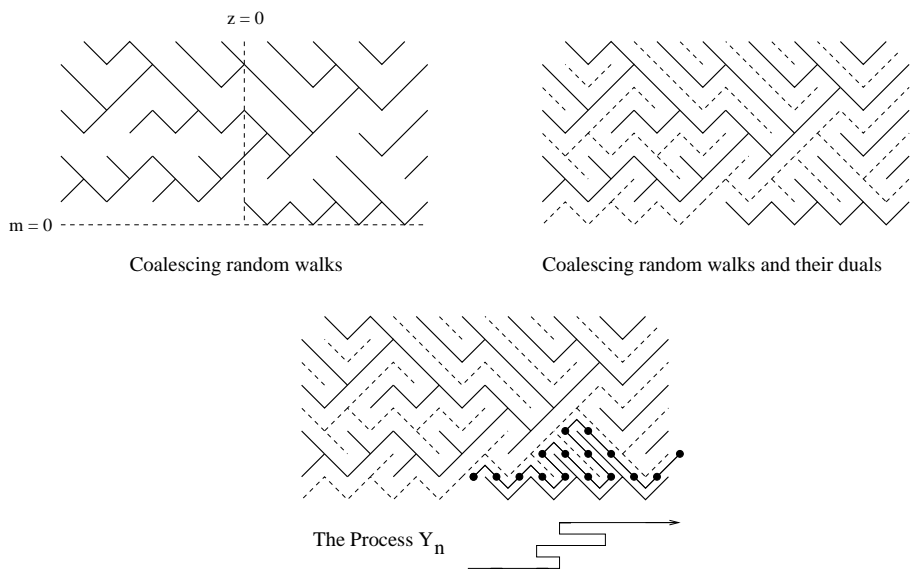

Coalescing random walks

Coalescing random walks and their duals

The Process $Y_n$

tion 11]. In particular, it may be described as an "infinitely negatively edge-reinforced random walk with initial occupation measure alternating between zero and one". To be more precise, give nearest neighbor edges of $z$ weight 1 if their center is at $\pm(1/2 + 2k)$ for $k = 0, 1, 2, \ldots$. Thus the two edges adjacent to zero are both labeled with a one, and, going away from zero in either direction, ones and zeros alternate. Now do a random walk that always chooses the less traveled edge, flipping a coin in the case of a tie (each crossing of an edge increases its weight by one).

The process $\{Y_n\}$ converges when rescaled to the process $\{X_t\}$ which is the scaling limit of the true self-avoiding walk. The limit operation in this case is more transparent: the coalescing simple random walks turn into coalescing Brownian motions. These Brownian motions are the local time processes given by the Ray-Knight theory. The construction of the process $\{X_t\}$ in [TW98] is in fact via these coalescing Brownian motions.

*The Stochastic Loewner Equation*

Suppose that the loop-erased random walk has a scaling limit. For specificity, it will be convenient to use the time reversal property of LERW and think of the walk as beginning on the boundary of a large disk and conditioned to hit the origin before returning to the boundary of the disk. The recursive $h$-process formulation (6.4) indicates that the infinitesimal future of such a limiting path would be a Brownian motion conditioned to avoid the path it has traced so far. Such conditioning, even if well defined, would seem to be complicated. But suppose, which is known about unconditioned Brownian motion and widely believed about many scaling limits, that the limiting LERW is conformally in-variant. The complement of the infinite past is simply connected, hence by the Riemann Mapping Theorem, it is conformally homeomorphic to the open unit disk with the present location mapping to a boundary point. The infinitesimal future in these coordinates is a Brownian motion conditioned immediately to enter the interior of the disk and stay there until it hits the origin. If we could compute in these coordinates, such conditioning would be routine.

In 2000, Schramm [Sch00] observed that such a conformal map may be computed via the classical **Löwner equation**. This is a differential equation satisfied by the conformal maps between a disk and the complement of a growing path inward from the boundary of the disk. More precisely, let $\beta$ be a compact simple path in the closed unit disk with one endpoint at zero and the other endpoint being the only point of $\beta$ on $\partial U$. Let $q : (-\infty, 0] \to \beta \setminus \{0\}$ be a parametrization of $\beta \setminus \{0\}$ and for each $t \leq 0$, let

$$f(t, z) : U \to U \setminus q([t, 0]) \tag{6.5}$$

be the unique conformal map fixing 0 and having positive real derivative at 0. Löwner[Löw23] proved that

**Theorem 6.14** (Löwner's Slit Mapping Theorem). *Given $\beta$, there is a parametrization $q$ and a continuous function $g : (-\infty, 0] \to \partial U$ such that the function $f : U \times (-\infty, 0] \to U$ in (6.5) satisfies the partial differential equation*

$$\frac{\partial f}{\partial t} = z \frac{g(t) + z}{g(t) - z} \frac{\partial f}{\partial z} \tag{6.6}$$

*with initial condition $f(z, 0) = z$.*

The point $q(t)$ is a boundary point of $U \setminus q([t, 0])$, so it corresponds under the Riemann map $f(t, \cdot)$ to a point on $\partial U$. It is easy to see this must be $g(t)$. Imagine that $\beta$ is the scaling limit of LERW started from the origin and stopped when it hits $\partial U$ (recurrence of two-dimensional random walk forces us to use a stopping construction). Since a Brownian motion conditioned to enter the interior of the disk has an angular component that is a simple Brownian motion, it is not too great a leap to believe that $g$ must be a Brownian motion on the circumference of $\partial U$, started from an arbitrary point, let us say 1. The solution to (6.6) exists for any $g$, that is, given $g$, we may recover the path $q$. We may then plug in for $g$ a Brownian motion with $\mathbb{E}B_t^2 = \kappa t$ for some scale parameter $\kappa$. We obtain what is known as the **radial $\mathsf{SLE}_\kappa$**.

More precisely, for any $\kappa > 0$, any simply connected open domain $D$, and any $x \in \partial D, y \in D$, there is a unique process $\mathsf{SLE}_\kappa(D; x, y)$ yielding a path $\beta$ as above from $x$ to $y$. We have constructed $\mathsf{SLE}_\kappa(D; 1, 0)$. This is sufficient because $\mathsf{SLE}_\kappa$ is invariant under conformal maps of the triple $(D; x, y)$. Letting $y$ approach $z \in \partial D$ gives a well defined limit known as **chordal $\mathsf{SLE}_\kappa(D; x, z)$**.

Lawler, Schramm and Werner have over a dozen substantial papers describing $\mathsf{SLE}_\kappa$ for various $\kappa$ and using $\mathsf{SLE}$ to analyze various scaling limits and solve some longstanding problems. A number of properties are proved in [RS05]. For example, $\mathsf{SLE}_\kappa$ is always a path, is self-avoiding if and only if $\kappa \leq 4$, and is space-filling when $\kappa \geq 8$. Regarding the question of whether SLE is the scaling limit of LERW, it was shown in [Sch00] that if LERW has a scaling limit and this is conformally invariant, then this limit is $SLE_2$. The conformally invariant limit was confirmed just a few years later:

**Theorem 6.15** ([LSW04, Theorem 1.3]). *Two-dimensional* LERW *stopped at the boundary of a disk has a scaling limit and this limit is conformally invariant. Consequently, the limit is* $SLE_2$.

In the same paper, Lawler, Schramm and Werner show that the peano curve separating an infinite uniform spanning tree from its dual has $SLE_8$ as its scaling limit. The $SLE_6$ is not self-avoiding, but its outer boundary is, up to an inessential transformation, the same as the outer boundary of a two-dimensional Brownian motion run until a certain stopping time. A recently announced result of Smirnov is that the interface between positive and negative clusters of the two-dimensional Ising model is an $SLE_3$. It is conjectured that the scaling limit of the classical self-avoiding random walk is $SLE_{8/3}$, the conjecture following if such a scaling limit can be proved to exist and be conformally invariant.

**Acknowledgements**

Thanks to M. Benaïm, G. Lawler, H. Mahmoud, S. Sheffield, B. Skyrms, B. Tóth and S. Volkov for comments on a preliminary draft. Some of the applications of urn models were collected by Tong Zhu in her masters thesis at the University of Pennsylvania [Zhu05].